\documentclass[10pt]{article}
\usepackage{amsmath}
\usepackage{amssymb}
\usepackage{mathrsfs}
\usepackage{amsfonts}
\usepackage{mathrsfs,amscd,amssymb,amsthm,amsmath,bm,graphicx}
\usepackage{graphicx}
\usepackage{cite}
\usepackage{color}
\makeatletter

\renewcommand{\@seccntformat}[1]{{\csname the#1\endcsname}{\normalsize.}\hspace{.5em}}
\makeatother

\def \[{\begin{equation}}
\def \]{\end{equation}}

\newtheorem{thm}{Theorem}[section]

\newtheorem{lem}[thm]{Lemma}

\newtheorem{cor}[thm]{Corollary}

\begin{document}
\setlength{\baselineskip}{13pt}
\begin{center}{\Large \bf  The normalized Laplacians and random walks of the parallel subdivision graphs
%\footnote{%The research is partially supported by National
%%Science Foundation of China (Grant No.  10671081)}
}

\vspace{4mm}

{\large Jing Zhao$^1$, Jia-Bao Liu$^{1,2,*}$, Ying-Ying Tan$^1$, Sakander Hayat$^3$}\vspace{2mm}

{\small $^1$School of Mathematics and Physics, Anhui Jianzhu
University, Hefei
230601, P.R. China\\
\small $^2$School of Mathematics, Southeast University, Nanjing 210096, P.R. China\\
$^3$ Faculty of Engineering Sciences, GIK Institute of Engineering Sciences and Technology, Topi 23460, Pakistan
}
%$^3$Department of Computer and Information Engineering, Wuhan Polytechnic University, Wuhan 430023, P.R. China.}
\vspace{2mm}
\end{center}

\footnotetext{E-mail address: zhaojing94823@163.com,
liujiabaoad@163.com,~tansusan1@ahjzu.edu.cn,~sakander1566@gmail.com.}

\footnotetext{* Corresponding author.}

 {\noindent{\bf Abstract.}\ \ The $k$-parallel subdivision graph $S_k(G)$ is generated from $G$ which each edge of $G$ is replaced by $k$ parallel paths of length 2. The $2k$-parallel subdivision graph $S_{2k}(G)$ is constructed from $G$ which each edge of $G$ is replaced by $k$ parallel paths of length 3. In this paper, the normalized Laplacian spectra of $S_k(G)$ and $S_{2k}(G)$ are given. They turn out that the multiplicities of the corresponding eigenvalues are only determined by $k$. As applications, the expected hitting time, the expected commute time and any two-points resistance distance between vertices $i$ and $j$ of $S_k(G)$, the normalized Laplacian spectra of $S_k(G)$ and $S_{2k}(G)$ with $r$ iterations are given. Moreover, the multiplicative degree Kirchhoff index, Kemeny's constant and the number of spanning tress of $S_k(G)$, $S_k^r(G)$, $S_{2k}(G)$ and $S_{2k}^r(G)$ are respectively obtained. Our results have generalized the previous works in Xie et al. and Guo et al. respectively.

\noindent{\bf Keywords}: Subdivision graph, Normalized Laplacian, Expected hitting time, Multiplicative degree Kirchhoff index, Kemeny's constant, Spanning tress

\noindent{\bf AMS subject classification:} 05C50,\ 60J20}

\section{Introduction}
\ \ \ \ Our convention in this paper is that the $(n,m)$-graph $G=(V(G),E(G))$ is the simple and connected graph with $n$ vertices and $m$ edges. For the associated notation of the graph theory, we follow \cite{Bondy}. One of the most critical fields in the graph theory is spectral graph theory. In particular, spectral graph theory is consist of algebraic spectral graph theory and analytic spectral graph theory. Spectral graph theory not only have pure mathematics properties, but also several vital applications in practical, examples include information-theoretic hashing of 3D objects\cite{T}, threshold selection in gene co-expression networks\cite{P}, pattern recognition, data mining and image matching\cite{J}. Thus, it is of great interests to determine the spectra of some graphs, especially for those large graphs which resulted from the small graph by graph operations in past two decades. For adjacency and Laplacian spectra of certain graph operations, see \cite{S1,S2,S3,D1,D2,D3,I1,Y1,G1,J1,X1,X2,C1,S4} and references therein. Alone this line, we consider the normalized Laplacian spectra of two parallel subdivision graphs in this paper, see Figure 1.

It is worthy to mentioned that $S_1(G)$ is the subdivision graph $S(G)$ if $k=1$ for $S_k(G)$ that considered in\cite{.P}. Besides, $S_2(G)$ is such diamond hierarchical graph which was constructed by Guo et al. \cite{Guo} if $k=2$ for $S_k(G)$. Moreover, while $k=1$ and $G$ is regular for $S_{2k}(G)$, then the graph $S_{2}(G)$ is the equivalent graph(generated by using clique-star transformation) of the stellated regular graphs\cite{Yang}.

In what follows, we will recall the definition of the normalized Laplacian. The formal definition is defined as
\begin{eqnarray*}
\mathscr{L}(G)=
\left\{
  \begin{array}{ll}
    \frac{-w(i,j)}{\sqrt{d(i)d(j)}}, & \hbox{$i$ is adjacent to $j$;} \\
    \frac{d(i)-w(i,i)}{d(i)}, & \hbox{$i=j, d(i)\neq0$;} \\
    ~~~~~~0, & \hbox{otherwise.}
  \end{array}
\right.
\end{eqnarray*}

\begin{figure}[htbp]
\centering\includegraphics[width=14cm,height=6cm]{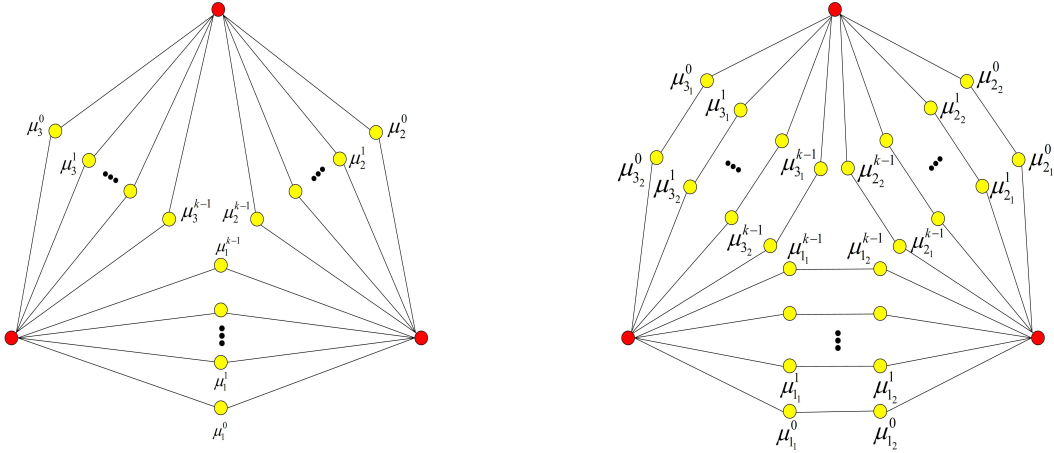}
\caption{(a) $S_{k}(G)$~~~~~~~~~~~~~~~~~~~~~~~~~~~~~~~~~~~~~~~~~~~~~(b) $S_{2k}(G)$}\label{1}
\end{figure}

In the definition defined as above, $w(i,j)$ is the weight of the edge $ij$ and $w(i,u)$ the vertex weight. In particular, $w(i,j)=1$ and $w(i,i)=0$ if vertex weight and edge weight are respectively 0 and 1. This yields
\begin{eqnarray}
\mathscr{L}(G)=
\left\{
  \begin{array}{ll}
    \frac{-1}{\sqrt{d(i)d(j)}}, & \hbox{$i$ is adjacent to $j$;} \\
    ~~~~~~1, & \hbox{$i=j, d(i)\neq0$;} \\
    ~~~~~~0, & \hbox{else.}
  \end{array}
\right.
\end{eqnarray}
Apparently, one has $\mathscr{L}(G)=D(G)^{-\frac{1}{2}}L(G)D(G)^{-\frac{1}{2}}$.

Recent years, Xie and Zhang \cite{P..,.P} respectively determined the normalized Laplacians of subdivisions and iterated triangulations graph. In 2017, the normalized Laplacians for quadrilateral graphs are considered by Li and Hou\cite{D.}. Pan et al.\cite{P.} computed the
normalized Laplacians of some graphs which involve graph transformations. For details, see \cite{J.,J..,Pan,Liu,He,Butter} and references therein.

A topic of a good deal of attention in network science is the notion of network criticality. some measures for network criticality are
depended on the paths in networks. Thus the shortest paths in networks are the focus of attention while we compute network criticality. But
the shortest paths are not considering the impacts of all the paths. The concepts of resistance distance and Kirchhoff index are such measures which involving all the information among paths in networks, one may refer to \cite{K}.

Resistance distance \cite{Klein} is defined by Klein and Randi\'c. Denote $\Omega_{ij}$ the resistance distance between vertices $i$ and $j$. In \cite{Ke}, Klein et al. defined the resistance distance sum rules, which is Kirchhoff index $Kf(G)$. The resistance distance based graph invariants have attracted some experts due to these applications in harmonic analysis, random walks, chemical and computer science. Considering the impacts of the vertex degree, Chen et al. \cite{Chen} defined the multiplicative degree Kirchhoff index $Kf^*(G)$ and found that there is a surprising relation between $Kf^*(G)$ and the eigenvalues of $\mathscr{L}(G)$.

The hitting time or first passage time $T_j$ of the vertex
$j$ is defined as the minimum steps of the random walk needs to reach that vertex. The expected hitting time or mean first passage time when the walk starts at $i$ is denoted by $E_iT_j(G)$. The expected commute time $C_{uv}(G)$ is defined as the expected time it takes for a particle to start from vertex $i$, reach vertex $j$, and then return to vertex $i$ in the graph $G$. The Kemeny's constant $Ke(G)$ is the expected number of steps needed for the transition from a starting vertex to a destination vertex, which is related to the finite ergodic Markov chains. Spanning trees are the subgraphs of $G$ which contain each edge in $G$, also an important quantity characterizing the reliability of a network, one may refer to \cite{Hunter,Yu}.

For a general graph, it is hard enough to determine the resistance distance between any points, unless the graph is very small or one knows the complete information of the graph. A. K. Chandra et al. in \cite{Chandra} built a relation between the expected hitting time and resistance distance of a graph. This provides a available method to calculate any two-points resistance distance, see\cite{Geo,exp,Law,Hit,ang}. The main results of this paper are presented as below.

\subsection{ The normalized Laplacian of $S_k(G)$ }
\ \ \ \ In this subsection, we will give the complete information for the normalized Laplacian of the $k$-parallel subdivision graph. Before proceeding, we shall disgress to define two functions as below.
\begin{eqnarray*}
f_1(x)=\frac{2+\sqrt{4-2x}}{2},~f_2(x)=\frac{2-\sqrt{4-2x}}{2}.
\end{eqnarray*}

\begin{thm} Assume that $G$ is a $(n,m)$-graph. Thus the eigenvalues of $\mathscr{L}(S_k(G))$ are as below.
\end{thm}
\begin{itemize}
\item $f_1(\lambda)$ and $f_2(\lambda)$ are eigenvalues of $\mathscr{L}(S_k(G))$, if $\lambda (\lambda\neq0)$ is an eigenvalue of $\mathscr{L}(G).$ Moreover, the multiplicities of $f_1(\lambda)$ and $f_2(\lambda)$ are same as $\lambda$.
\item $0$ and $2$ are the eigenvalues of $\mathscr{L}(S_k(G))$ with the multiplicity $1$.
\item $\mathscr{L}(S_k(G))$ has the eigenvalue 1. Its multiplicity is $km-n+2$, if $G$ is bipartite and $km-n$ else.
\end{itemize}

\subsection{ The normalized Laplacian of $S_{2k}(G)$ }
\ \ \ \ In this subsection, the normalized Laplacian eigenvalues of the $2k$-parallel subdivision graph are proposed. In addition, let $\mu_1,$ $\mu_2$ and $\mu_3$ be the roots of the equation $4\mu^3-12\mu^2+9\mu-\lambda=0$.

\begin{thm} Assume that $G$ is a $(n,m)$-graph. Therefore the eigenvalues of $\mathscr{L}(S_{2k}(G))$ are as follows.
\end{thm}
\begin{itemize}
\item $\mu_1,$ $\mu_2$ and $\mu_3$ are eigenvalues of $\mathscr{L}(S_{2k}(G))$, if $\lambda (\lambda\neq0,2)$ is an eigenvalue of $\mathscr{L}(G).$ The multiplicities of $\mu_1,$ $\mu_2$ and $\mu_3$ are same as $\lambda$.
\item $\mathscr{L}(S_{2k}(G))$ has an eigenvalue $0$ with the multiplicity $1$. $\mathscr{L}(S_{2k}(G))$ has an eigenvalue $2$ and its multiplicity is $1$, if $G$ is bipartite.
\item $\mathscr{L}(S_{2k}(G))$ has the eigenvalues $\frac{1}{2}$ and $\frac{3}{2}$. Their multiplicities are respectively $km-n$ and $km-n+2$, if $G$ is non-bipartite.
\item $\mathscr{L}(S_{2k}(G))$ has the eigenvalues $\frac{1}{2}$ and $\frac{3}{2}$. Their multiplicities are both $km-n+2$, if $G$ is bipartite.
\end{itemize}

The rest parts are summarized as below. We propose some lemmas in section 2. The proofs of the main results are provided in section 3. The spectra of $\mathscr{L}(S_{k}(G))$ and $\mathscr{L}(S_{2k}(G))$ are determined in sections 4 and 5. As byproduct, the formulas for $E_iT_j(S_{k}(G))$, $r_{ij}(S_{k}(G))$, $Kf^*(S_{k}(G))$($Kf^*(S_{2k}(G))$ resp.), $Ke(S_{k}(G))$($Kf^*(S_{2k}(G))$ resp.) and $\tau(S_{k}(G))$($\tau(S_{2k}(G))$ resp.) are respectively given.

\section{ Preliminaries}
In this section, we put some lemmas that used in applications section. Let $0=\lambda_1<\lambda_2\leq\cdots\leq\lambda_n$ be the normalized Laplacian eigenvalues of $\mathscr{L}(G)$.

In 2007, Chen and Zhang found that $Kf^*(G)$ is closely related to the eigenvalues of $\mathscr{L}(S_k(G))$ in $G$.
\begin{lem} \cite{Chen} Assume that $G$ is a $(n,m)$-graph. Then
$Kf^*(G)=2m\sum_{i=2}^n\frac{1}{\lambda_i}.$
\end{lem}

Kemeny and Snell in 1960 introduced the graph invariant which is called Kemeny's constant. It also related to the eigenvalues of $\mathscr{L}(S_k(G))$.
\begin{lem} \cite{Butter} Assume that $G$ is a $(n,m)$-graph. Then
$Ke(G)=\sum_{i=2}^n\frac{1}{\lambda_i}.$
\end{lem}

Combining Lemmas 2.1 and 2.2, one finds $Kf^*(G)=2m\cdot Ke(G)$. Let $B(G)$ denote the vertex-edge incident matrix and $r(B(G))$ the rank of $B(G)$, one has the following.
\begin{lem} \cite{Butter} Assume that $G$ is a $(n,m)$-graph. Thus
\begin{eqnarray*}
r(B(G))=
\left\{
  \begin{array}{ll}
    n, & \hbox{if $G$ is non-bipartite;} \\
    n-1, & \hbox{if $G$ is bipartite.}
  \end{array}
\right.
\end{eqnarray*}
\end{lem}

\begin{lem} \cite{F.R} Assume that $G$ is a $(n,m)$-graph. One obtains
\begin{eqnarray*}
\prod_{i=1}^nd(v_i)\prod_{i=2}^n\lambda_i=2m\tau(G),
\end{eqnarray*}
where
$\tau(G)$ is the number of spanning trees of $G$.
\end{lem}

At this point, one considers the normalized adjacency matrix $N(G)$, namely
\begin{eqnarray*}
N(G)=D(G)^{-\frac{1}{2}}A(G)D(G)^{-\frac{1}{2}}=D(G)^{\frac{1}{2}}(I-\mathscr{L}(G))D(G)^{-\frac{1}{2}},
\end{eqnarray*}
where $I$ is the identity matrix.

The above equation leads that matrix $N(G)$ has the same eigenvalues as matrix $I-\mathscr{L}(G)$. Assume that the eigenvalues of $I-\mathscr{L}(G)$ has the orthogonal eigenvectors $v_1,v_2,\ldots,v_n$. Namely
$$v_i=(v_{i1},v_{i2},\ldots,v_{in})^T,~i=1,2,\ldots,n.$$

One knows that $E=(v_1,v_2,\ldots,v_n)$ is the orthogonal matrix, then
\begin{eqnarray*}
\sum_{k=1}^{n}v_{ik}v_{jk}=\sum_{k=1}^{n}v_{ki}v_{kj}=\left\{
                                                          \begin{array}{ll}
                                                            1, & \hbox{if $i=j$;} \\
                                                            0, & \hbox{else.}
                                                          \end{array}
                                                        \right.
\end{eqnarray*}

In particular,
\begin{eqnarray*}
v_1=\bigg(\sqrt{d_1/2m}, \sqrt{d_2/2m}, ... , \sqrt{d_n/2m}\bigg).
\end{eqnarray*}
Let $G$ be a bipartite graph with a vertex partition $V(G)=V_1\bigcup V_2$. Thus
\begin{eqnarray*}
v_{ni}=\sqrt{d_i/2m}, i\in V_1; v_{nj}=-\sqrt{d_j/2m}, j\in V_2.
\end{eqnarray*}

\begin{lem} \cite{Ran} Assume that $G$ is a $(n,m)$-graph. One obtains
\begin{eqnarray*}
E_iT_j(G)=2m\sum_{a=2}^{n}\frac{1}{1-\lambda_a}
\bigg(\frac{v_{aj}^2}{d_j}+\frac{v_{ai}v_{aj}}{\sqrt{d_id_j}}\bigg) .
\end{eqnarray*}
\end{lem}

\begin{lem} \cite{Chandra} Assume that $G$ is a $(n,m)$-graph. One gets
\begin{eqnarray*}
E_iT_j(G)+E_jT_i(G)=2m\Omega_{ij}(G).
\end{eqnarray*}
\end{lem}

The following lemma built a connection between the expected hitting time and expected commute time.
\begin{lem} \cite{Chandra} Assume that $G$ is a $(n,m)$-graph. One gets
\begin{eqnarray*}
E_iT_j(G)+E_jT_i(G)=C_{uv}(G).
\end{eqnarray*}
\end{lem}
Certainly, one reaches $C_{uv}(G)=2m\Omega_{ij}(G)$.

\section{ The proofs of Theorems 1.1 and 1.2}
\ \ \ \ The most focuses in this section are on determining the eigenvalues of $\mathscr{L}(S_{k}(G))$ and $\mathscr{L}(S_{2k}(G))$. Before we begin, we turn to prove two lemmas as basic but necessary in the proofs.

\begin{lem} $\mathscr{L}(G)$ has an eigenvalue $4\sigma-2\sigma^2$and its multiplicity is same as $\sigma$, if $\sigma(\sigma\neq1)$ is an eigenvalue of $\mathscr{L}(S_{k}(G))$.
\end{lem}
\noindent{\bf Proof.} According to the construction of $S_{k}(G)$, let $N_0=\{\mu_1^0,\mu_2^0,\ldots,\mu_m^0\}$, $N_1=\{\mu_1^1,\mu_2^1,\ldots,\mu_m^1\}$, \ldots, $N_{k-1}=\{\mu_1^{k-1},\mu_2^{k-1},\ldots,\mu_m^{k-1}\}.$ Evidently, $N=N_0\bigcup N_1\bigcup \cdots\bigcup N_{k-1}$. At this place, we recall that
\begin{eqnarray}
d_w(S_{k}(G))=
\begin{cases}
k\cdot d_w(G),       & w\in V(G); \\
2,           & w\in N.
\end{cases}
\end{eqnarray}

Set $\varepsilon=\{\varepsilon_1,\varepsilon_2,\ldots,\varepsilon_{n+km}\}$ as an eigenvector corresponding to the eigenvalue $\sigma$
of $S_{k}(G)$. Together the definition of normalized Laplacian with Eq.(1.1), one gets the following equation for any vertex $w$ in
$S_{k}(G)$.
\begin{eqnarray}
(1-\sigma)\varepsilon_w=\sum_{u\thicksim w}\frac{1}{\sqrt{d_w(S_{k})d_u(S_{k})}}\varepsilon_u,
\end{eqnarray}
where denote $S_{k}(G)$ by $S_{k}$ for simplify.

Suppose that $N_l(u)$ is the set which consists of all neighbors of $u$ in $N_l~(l=0,1,\ldots,k-1)$. The set of all neighbors of $u$ in $G$ is denoted by $N_G(u)$.
For any vertex $u$ in $G$ and Eqs.(3.2) and (3.3), one has
\begin{eqnarray}
(1-\sigma)\varepsilon_u=\sum_{l=0}^{k-1}\sum_{u_i^l\in N_l(u)}\frac{1}{\sqrt{d_{u_i^l}(S_{k})d_u(S_{k})}}\varepsilon_{u_i^l}
=\sum_{l=0}^{k-1}\sum_{u_i^l\in N_l(u)}\frac{1}{\sqrt{2k\cdot d_u(G)}}\varepsilon_{u_i^l}.
\end{eqnarray}

For any vertex $u_i^l\in N_l(u)$, one obtains
\begin{eqnarray}
(1-\sigma)\varepsilon_{u_i^l}=\frac{1}{\sqrt{d_u(S_{k})d_{u_i^l}(S_{k})}}\varepsilon_{u}+
\frac{1}{\sqrt{d_v(S_{k})d_{u_i^l}(S_{k})}}\varepsilon_{v}
=\frac{1}{\sqrt{2k\cdot d_u(G)}}\varepsilon_{u}+
\frac{1}{\sqrt{2k\cdot d_v(G)}}\varepsilon_{v},
\end{eqnarray}
where vertex $v$ is adjacent to vertex $u$ in $G$.

Combining Eq.(3.4) and Eq.(3.5), we have
\begin{eqnarray}
(1-\sigma)^2\varepsilon_{u}&=&\sum_{l=0}^{k-1}\sum_{u_i^l\in N_l(u)}\bigg(\frac{1}{2k\cdot d_u(G)}\varepsilon_{u}+
\frac{1}{2k\sqrt{d_u(G)d_v(G)}}\varepsilon_{v}\bigg) \nonumber \\
&=&k\cdot\sum_{v\in N_G(u)}\bigg(\frac{1}{2k\cdot d_u(G)}\varepsilon_{u}+
\frac{1}{2k\sqrt{d_u(G)d_v(G)}}\varepsilon_{v}\bigg) \nonumber \\
&=&\frac{1}{2}\varepsilon_{u}+\frac{1}{2}\sum_{v\in N_G(u)}\frac{1}{\sqrt{d_u(G)d_v(G)}}\varepsilon_{v},
\end{eqnarray}
for $\sigma\neq1$.

By a straightforward calculation of Eq.(3.6), it gives
\begin{eqnarray*}
\big[1-(4\sigma-2\sigma^2)\big]\varepsilon_{u}=\sum_{v\in N_G(u)}\frac{1}{\sqrt{d_u(G)d_v(G)}}\varepsilon_{v},
\end{eqnarray*}
for $\sigma\neq1$.

Based on Eq.(3.3), one knows $4\sigma-2\sigma^2$ is an eigenvalue of $\mathscr{L}(G)$. Furthermore, there exits a bijection between $\sigma$
and $4\sigma-2\sigma^2$, thus the multiplicity of $4\sigma-2\sigma^2$ is same as that of $\sigma$.\\
This has completed the proof.\hfill\rule{1ex}{1ex}

\begin{lem} $\mathscr{L}(G)$ has an eigenvalue $\zeta(4\zeta^2-12\zeta+9)$ and its multiplicity is same as $\zeta$, if $\zeta(\zeta\neq\frac{1}{2},\frac{3}{2})$ is an eigenvalue of $\mathscr{L}(S_{2k}(G))$.
\end{lem}
\noindent{\bf Proof.} For the graph $S_{2k}(G)$, let $M_0=\{\mu_{1_1}^0,\mu_{1_2}^0,\ldots,\mu_{m_1}^0,\mu_{m_2}^0\}$, $N_1=\{\mu_{1_1}^1,\mu_{1_2}^1,\ldots,\mu_{m_1}^1,\mu_{m_2}^1\}$, \ldots, $N_{k-1}=\{\mu_{1_1}^{k-1},\mu_{1_2}^{k-1},\ldots,\mu_{m_1}^{k-1},\mu_{m_2}^{k-1}\}.$ Obviously, $M=M_0\bigcup M_1\bigcup \cdots\bigcup M_{k-1}$. At this point, we review that
\begin{eqnarray}
d_w(S_{2k}(G))=
\begin{cases}
k\cdot d_w(G),       & w\in V(G); \\
2,           & w\in M.
\end{cases}
\end{eqnarray}

Set $\xi=\{\xi_1,\xi_2,\ldots,\xi_{n+2km}\}$ as an eigenvector corresponding to the eigenvalue $\zeta$
of $S_{2k}(G)$. Together the definition of normalized Laplacian with Eq.(1.1), one gets the following equation for any vertex $s$ in
$S_{2k}(G)$.
\begin{eqnarray}
(1-\zeta)\xi_s=\sum_{v\thicksim s}\frac{1}{\sqrt{d_s(S_{2k})d_v(S_{2k})}}\xi_v,
\end{eqnarray}
where denote $S_{2k}(G)$ by $S_{2k}$ for short.

Assume that $M_l(u)$ is the set which consists of all neighbors of $u$ in $M_l~(l=0,1,\ldots,k-1)$. Denote the set of all neighbors of $u$ in $G$ by $N_G(u)$.
Any vertex $u$ in $G$ and Eq.(3.7) imply that
\begin{eqnarray}
(1-\zeta)\xi_u=\sum_{l=0}^{k-1}\sum_{u_{i_1}^l\in M_l(u)}\frac{1}{\sqrt{d_{u_{i_1}^l}(S_{2k})d_u(S_{2k})}}\xi_{u_{i_1}^l}
=\sum_{l=0}^{k-1}\sum_{u_{i_1}^l\in M_l(u)}\frac{1}{\sqrt{2k\cdot d_u(G)}}\xi_{u_{i_1}^l}.
\end{eqnarray}

For any $u_{i_1}^l\in M_l(u)$, we obtain
\begin{eqnarray}
(1-\zeta)\xi_{u_{i_1}^l}&=&\frac{1}{\sqrt{d_u(S_{2k})d_{u_{i_1}^l}(S_{2k})}}\xi_u+\frac{1}{\sqrt{d_{u_{i_1}^l}(S_{2k})d_{u_{i_2}^l}
(S_{2k})}}\xi_{u_{i_2}^l}\nonumber\\
&=&\frac{1}{\sqrt{2k\cdot d_u(G)}}\xi_u+\frac{1}{2}\xi_{u_{i_2}^l},
\end{eqnarray}
where vertex $u_{i_2}^l$ in $M_l$ is adjacent to vertex $u_{i_1}^l$.

Similarly, for any $u_{i_2}^l\in M_l$, this yields
\begin{eqnarray}
(1-\zeta)\xi_{u_{i_2}^l}&=&\frac{1}{\sqrt{d_{u_{i_1}^l}(S_{2k})d_{u_{i_2}^l}(S_{2k})}}\xi_{u_{i_1}^l}+
\frac{1}{\sqrt{d_v(S_{2k})d_{u_{i_2}^l}(S_{2k})}}\xi_v\nonumber\\
&=&\frac{1}{2}\xi_{u_{i_1}^l}+\frac{1}{\sqrt{2k\cdot d_v(G)}}\xi_v.
\end{eqnarray}

Combining Eq.(3.11) and Eq.(3.10), one arrives at
\begin{eqnarray}
2\big(\zeta-\frac{1}{2}\big)\big(\zeta-\frac{3}{2}\big)\xi_{u_{i_1}^l}=\frac{2(1-\zeta)}{\sqrt{2k\cdot d_u(G)}}\xi_u
+\frac{1}{\sqrt{2k\cdot d_v(G)}}\xi_v.
\end{eqnarray}

Substituting Eq.(3.12) to Eq.(3.9), one gets
\begin{eqnarray}
2\big(\zeta-\frac{1}{2}\big)\big(\zeta-\frac{3}{2}\big)\big(1-\zeta\big)\xi_{u}&=&k\cdot\sum_{v\in N_G(u)}\bigg(\frac{1-\zeta}{k\cdot d_u(G)}\xi_{u}+\frac{1}{2k\sqrt{d_u(G)d_v(G)}}\xi_{v}\bigg)\nonumber\\
&=&(1-\zeta)\xi_u+\frac{1}{2}\sum_{v\in N_G(u)}\frac{1}{\sqrt{d_u(G)d_v(G)}}\xi_{v},
\end{eqnarray}
where $\zeta\neq\frac{1}{2},\frac{3}{2}$.

By a explicit analysis of Eq.(3.13), that is
\begin{eqnarray*}
\big[1-\zeta(4\zeta^2-12\zeta+9)\big]\xi_{u}=\sum_{v\in N_G(u)}\frac{1}{\sqrt{d_u(G)d_v(G)}}\xi_{v},
\end{eqnarray*}
for $\zeta\neq\frac{1}{2},\frac{3}{2}$.

According to Eq.(3.8), one gets $\zeta(4\zeta^2-12\zeta+9)$ is an eigenvalue of $\mathscr{L}(G)$. Moreover, there exits a bijection between $\mu$ and $\zeta\neq\frac{1}{2},\frac{3}{2}$, thus the multiplicity of $\zeta\neq\frac{1}{2},\frac{3}{2}$ is same as that of $\zeta$.\\
This completes the proof.\hfill\rule{1ex}{1ex}

Now, we proceed by going into more details on the proofs of the normalized Laplacians of $S_{k}(G)$ and $S_{2k}(G)$. Our first goal in the following is determined the eigenvalues of $\mathscr{L}(S_{k}(G))$.

\subsection{ The proof of Theorem 1.1 }
\ \ \ \ At this point, we slightly to observe those two functions that defined in the subsection 1.2 as follows.
\begin{eqnarray}
f_1(x)=\frac{2+\sqrt{4-2x}}{2},~f_2(x)=\frac{2-\sqrt{4-2x}}{2}.
\end{eqnarray}

 Assume that $\sigma(\sigma\neq1)$ is an eigenvalue of $\mathscr{L}(S_{k}(G))$ and $\lambda$ an eigenvalue of $\mathscr{L}(G)$, then by Lemma 3.1, one has $\lambda=4\sigma-2\sigma^2$. Put it in another way, one obtains $\sigma=\frac{2\pm\sqrt{4-2\lambda}}{2},~\sigma\neq1$.

It is routine to check that $\mathscr{L}(G)$ has an eigenvalue $0$ and its multiplicity is $1$. Then substituting $0$ into Eq.(3.14), we get $\sigma=0, 2$ with the multiplicity $1$.

The graph $S_{k}(G)$ is bipartite whether $G$ is bipartite or not. By Lemma 3.1, one obtains the rest of the eigenvalues of $S_{k}(G)$ are $1$. Further, $1$ is with the multiplicity $km-n+2$ if $G$ is bipartite and $km-n$ otherwise.

\subsection{ The proof of Theorem 1.2 }
\ \ \ \
 Assume that $\mathscr{L}(S_{2k}(G))$ has an eigenvalue $\zeta(\zeta\neq\frac{1}{2},\frac{3}{2})$ and $\mathscr{L}(G)$ has an eigenvalue $\lambda$, then by Lemma 3.2, one has $\lambda=\zeta(4\zeta^2-12\zeta+9)$. Namely, one gets
\begin{eqnarray}
    4\zeta^3-12\zeta^2+9\zeta-\lambda=0,~\zeta\neq\frac{1}{2},\frac{3}{2}.
\end{eqnarray}

 It is worthy to mention that $0$ is a normalized Laplacian eigenvalue of $G$ with the multiplicity $1$. Combining Eq.(3.15) and $\lambda=0$, this gives $\zeta=0,\frac{3}{2},\frac{3}{2}$.

 Note that $\mathscr{L}(G)$ has an eigenvalue $2$ with the multiplicity $1$, when $G$ is bipartite. Substituting $\lambda=2$ into Eq.(3.15), then one arrives at $\zeta=2,\frac{1}{2},\frac{1}{2}$.

 Assume that $G$ is non-bipartite. Combining Eq.(3.12) and $\zeta=\frac{1}{2}$, one gets
\begin{eqnarray}
    \frac{\xi_u}{\sqrt{d_u(G)}}=-\frac{\xi_v}{\sqrt{d_v(G)}}.
\end{eqnarray}

There must be an odd cycle in $G$ due to $G$ is non-bipartite. Assume that the length of that odd cycle is $p$ and those vertices are named $t_1,t_2,\ldots,t_p$. Thus
 \begin{eqnarray*}
    \frac{\xi_{t_1}}{\sqrt{d_{t_1}(G)}}=-\frac{\xi_{t_2}}{\sqrt{d_{t_2}(G)}}=\cdots=\frac{\xi_{t_p}}{\sqrt{d_{t_p}(G)}}
    =-\frac{\xi_{t_1}}{\sqrt{d_{t_1}(G)}}.
\end{eqnarray*}

From the above equation, one has $\xi_u=0,~u\in V(G)$. Alone with Eq.(3.9), we have
\begin{eqnarray*}
\sum_{l=0}^{k-1}\sum_{u_{i_1}^l\in M_l(u)}\frac{1}{\sqrt{2k\cdot d_u(G)}}\xi_{u_{i_1}^l}=0,
\end{eqnarray*}
namely,
\begin{eqnarray}
\sum_{q\in M(u) }\xi_{q}=0,~M(u)=\bigcup_{l=0}^{k-1}M_l(u).
\end{eqnarray}

Substituting $\zeta=\frac{1}{2}$ and $\xi_u=0,~u\in V(G)$ into Eq.(3.10), one obtains $\xi_{u_{i_1}^l}-\xi_{u_{i_2}^l}=0,~l=0,1,\ldots,k-1$.
$$\xi_{u_{1_1}^0}=\xi_{u_{1_2}^0}=y_1^0,~\xi_{u_{1_1}^1}=\xi_{u_{1_2}^1}=y_1^1,\ldots,\xi_{u_{1_1}^{k-1}}=\xi_{u_{1_2}^{k-1}}=y_1^{k-1},$$
$$\xi_{u_{2_1}^0}=\xi_{u_{2_2}^0}=y_2^0,~\xi_{u_{2_1}^1}=\xi_{u_{2_2}^1}=y_2^1,\ldots,\xi_{u_{2_1}^{k-1}}=\xi_{u_{2_2}^{k-1}}=y_2^{k-1},$$
$$\vdots$$
$$\xi_{u_{m_1}^0}=\xi_{u_{m_2}^0}=y_m^0,~\xi_{u_{m_1}^1}=\xi_{u_{m_2}^1}=y_m^1,\ldots,\xi_{u_{m_1}^{k-1}}=\xi_{u_{m_2}^{k-1}}=y_m^{k-1}.$$

Let
\begin{eqnarray*}
y_0=(y_1^0,y_2^0,\ldots,y_m^0)^T,~y_1=(y_1^1,y_2^1,\ldots,y_m^1)^T, \ldots,~y_{k-1}=(y_1^{k-1},y_2^{k-1},\ldots,y_m^{k-1})^T.
\end{eqnarray*}
Assume that $B(G)=(\alpha_1, \alpha_2,\ldots, \alpha_n)^T$ is the incident matrix of $G$ and $y=(y_0,y_1, \ldots, y_{k-1})^T$.
According to Eq.(3.17), one can write it by another way, that is
\begin{eqnarray}
Cy=0,
\end{eqnarray}
where $C=(\beta_1, \beta_2,\ldots, \beta_n)^T$, $\beta_1=(\underbrace{\alpha_1, \alpha_1,\ldots, \alpha_1}_k)$,$\beta_2=(\underbrace{\alpha_2, \alpha_2,\ldots, \alpha_2}_k)$,\ldots,$\beta_n=(\underbrace{\alpha_n, \alpha_n,\ldots, \alpha_n}_k)$.

Evidently, $r(B(G))=r(C)$. Thus, Eq.(3.18) has $km-n$ linearly independent solutions based on Lemma (2.3). Put it in another way, the multiplicity of the eigenvalue $\frac{1}{2}$ of $\mathscr{L}(S_{2k}(G))$ is $km-n$. Therefore, the multiplicity of the eigenvalue $\frac{3}{2}$ of $\mathscr{L}(S_{2k}(G))$ is $km-n+2$.

The graph $S_{2k}(G)$ is bipartite when $G$ is bipartite, then . As we know, the eigenvalues of $\mathscr{L}(S_{2k}(G))$ are symmetric about $1$, thus the multiplicities of $\frac{1}{2}$ and $\frac{3}{2}$ are equal to $km-n+2$.

\section{Applications of Theorem 1.1}
\ \ \ \ In what follows, the expected hitting time and any two-points resistance distance between any vertices $i$ and $j$ of $S_{k}(G)$ are given. At here, we disgress to propose the adjacency and degree diagonal matrices of $S_{k}(G)$.

\begin{eqnarray*}
A(S_{k}(G))=\left(
                \begin{array}{cccc}
                  0 & B(G) & \cdots & B(G) \\
                  B(G)^T & 0 & \cdots & 0 \\
                  \vdots & \vdots & \ddots & \vdots \\
                  B(G)^T & 0 & \cdots & 0 \\
                \end{array}
              \right),
D(S_{k}^r(G))=\left(
                \begin{array}{cccc}
                  kD(G) & 0 & \cdots & 0 \\
                  0 & 2I_m & \cdots & 0 \\
                  \vdots & \vdots & \ddots & \vdots \\
                  0 & 0 & \cdots & 2I_m \\
                \end{array}
              \right).
\end{eqnarray*}

One arrives at
\begin{eqnarray*}
N(S_{k}(G))&=&D(G)^{-\frac{1}{2}}A(G)D(G)^{-\frac{1}{2}}\\&=&\left(
                                                  \begin{array}{cccc}
                                                    0 & \frac{1}{\sqrt{2k}}D(G)^{-\frac{1}{2}}B(G) & \cdots & \frac{1}{\sqrt{2k}}D(G)^{-\frac{1}{2}}B(G) \\
                                                    \frac{1}{\sqrt{2k}}B(G)^TD(G)^{-\frac{1}{2}} & 0 & \cdots & 0 \\
                                                    \vdots & \vdots & \ddots & \vdots \\
                                                    \frac{1}{\sqrt{2k}}B(G)^TD(G)^{-\frac{1}{2}} & 0 & \cdots & 0 \\
                                                  \end{array}
                                                \right).
\end{eqnarray*}

In the following lemma, we proposed the orthonormal eigenvectors of the corresponding eigenvalues in terms of Theorem 1.1 of $N(S_{k}(G))$.
\begin{lem}
Assume that $G$ is a $(n,m)$-graph. One has
\begin{itemize}
  \item $N(S_{k}(G))$ has the eigenvalues $\pm\sqrt{\frac{1+\lambda_a}{2}},~a=1,2,\ldots,n$ and 0 with multiplicities respectively 1 and $km-n$, if $G$ is non-bipartite graph. Then the corresponding orthonormal eigenvectors are
      \begin{eqnarray*}
\frac{1}{\sqrt{2}}\left(
                    \begin{array}{c}
                      v_i \\
                    \pm\frac{1}{\sqrt{k(1+\lambda_i)}}B^TD^{-\frac{1}{2}}v_i \\
                    \vdots \\
                    \pm\frac{1}{\sqrt{k(1+\lambda_i)}}B^TD^{-\frac{1}{2}}v_i
                    \end{array}
                  \right),~i=1,2,\ldots,n;
\end{eqnarray*}
    \begin{eqnarray*}
\left(
                    \begin{array}{c}
                      0 \\
                    \frac{1}{\sqrt{2(1+\lambda_i)}}B^TD^{-\frac{1}{2}}v_i \\
                    -\frac{1}{\sqrt{2(1+\lambda_i)}}B^TD^{-\frac{1}{2}}v_i \\
                    0\\
                    0 \\
                    \vdots \\
                    0 \\
                    \end{array}
                  \right),
\left(
                    \begin{array}{c}
                      0 \\
                    \frac{1}{\sqrt{6(1+\lambda_i)}}B^TD^{-\frac{1}{2}}v_i \\
                    \frac{1}{\sqrt{6(1+\lambda_i)}}B^TD^{-\frac{1}{2}}v_i \\
                    -\sqrt{\frac{2}{3(1+\lambda_i)}}B^TD^{-\frac{1}{2}}v_i\\
                    0 \\
                    \vdots \\
                    0 \\
                    \end{array}
                  \right),\ldots,
\left(
                    \begin{array}{c}
                      0 \\
                    \sqrt{\frac{1}{k(k-1)(1+\lambda_i)}}B^TD^{-\frac{1}{2}}v_i \\
                    \sqrt{\frac{1}{k(k-1)(1+\lambda_i)}}B^TD^{-\frac{1}{2}}v_i \\
                    \sqrt{\frac{1}{k(k-1)(1+\lambda_i)}}B^TD^{-\frac{1}{2}}v_i\\
                    \sqrt{\frac{1}{k(k-1)(1+\lambda_i)}}B^TD^{-\frac{1}{2}}v_i \\
                    \vdots \\
                    -\sqrt{\frac{k-1}{k(1+\lambda_i)}}B^TD^{-\frac{1}{2}}v_i \\
                    \end{array}
                  \right);
\end{eqnarray*}
   \begin{eqnarray*}
\left(
                    \begin{array}{c}
                      0 \\
                    y_z \\
                    0 \\
                    \vdots\\
                    0 \\
                    \end{array}
                  \right),
\left(
                    \begin{array}{c}
                      0 \\
                    0 \\
                    y_z \\
                    \vdots\\
                    0 \\
                    \end{array}
                  \right),\ldots,
\left(
                    \begin{array}{c}
                      0 \\
                    0 \\
                    0 \\
                    \vdots\\
                    y_z \\
                    \end{array}
                  \right),~z=1,2,\ldots,m-n,
\end{eqnarray*}
where $(y_1,y_2,\ldots,y_{m-n})$ is an orthonormal basis of the kernel space of the matrix $B(G)$.
  \item $N(S_{k}(G))$ has the eigenvalues $\pm\sqrt{\frac{1+\lambda_a}{2}},~a=1,2,\ldots,n-1$ and 0 with multiplicities respectively 1 and $km-n+2$, if $G$ is bipartite graph. Then the corresponding orthonormal eigenvectors are
      \begin{eqnarray*}
\frac{1}{\sqrt{2}}\left(
                    \begin{array}{c}
                      v_i \\
                    \pm\frac{1}{\sqrt{k(1+\lambda_i)}}B^TD^{-\frac{1}{2}}v_i \\
                    \vdots \\
                    \pm\frac{1}{\sqrt{k(1+\lambda_i)}}B^TD^{-\frac{1}{2}}v_i
                    \end{array}
                  \right),~i=1,2,\ldots,n-1;
\end{eqnarray*}
    \begin{eqnarray*}
\left(
                    \begin{array}{c}
                      0 \\
                    \frac{1}{\sqrt{2(1+\lambda_i)}}B^TD^{-\frac{1}{2}}v_i \\
                    -\frac{1}{\sqrt{2(1+\lambda_i)}}B^TD^{-\frac{1}{2}}v_i \\
                    0\\
                    0 \\
                    \vdots \\
                    0 \\
                    \end{array}
                  \right),
\left(
                    \begin{array}{c}
                      0 \\
                    \frac{1}{\sqrt{6(1+\lambda_i)}}B^TD^{-\frac{1}{2}}v_i \\
                    \frac{1}{\sqrt{6(1+\lambda_i)}}B^TD^{-\frac{1}{2}}v_i \\
                    -\sqrt{\frac{2}{3(1+\lambda_i)}}B^TD^{-\frac{1}{2}}v_i\\
                    0 \\
                    \vdots \\
                    0 \\
                    \end{array}
                  \right),\ldots,
\left(
                    \begin{array}{c}
                      0 \\
                    \sqrt{\frac{1}{k(k-1)(1+\lambda_i)}}B^TD^{-\frac{1}{2}}v_i \\
                    \sqrt{\frac{1}{k(k-1)(1+\lambda_i)}}B^TD^{-\frac{1}{2}}v_i \\
                    \sqrt{\frac{1}{k(k-1)(1+\lambda_i)}}B^TD^{-\frac{1}{2}}v_i\\
                    \sqrt{\frac{1}{k(k-1)(1+\lambda_i)}}B^TD^{-\frac{1}{2}}v_i \\
                    \vdots \\
                    -\sqrt{\frac{k-1}{k(1+\lambda_i)}}B^TD^{-\frac{1}{2}}v_i \\
                    \end{array}
                  \right);
\end{eqnarray*}
   \begin{eqnarray*}
   \left(
                    \begin{array}{c}
                      v_n \\
                    0 \\
                    0 \\
                    \vdots\\
                    0 \\
                    \end{array}
                  \right),
\left(
                    \begin{array}{c}
                      0 \\
                    y_z \\
                    0 \\
                    \vdots\\
                    0 \\
                    \end{array}
                  \right),
\left(
                    \begin{array}{c}
                      0 \\
                    0 \\
                    y_z \\
                    \vdots\\
                    0 \\
                    \end{array}
                  \right),\ldots,
\left(
                    \begin{array}{c}
                      0 \\
                    0 \\
                    0 \\
                    \vdots\\
                    y_z \\
                    \end{array}
                  \right),~z=1,2,\ldots,m-n+1,
\end{eqnarray*}
where $(y_1,y_2,\ldots,y_{m-n+1})$ is an orthonormal basis of the kernel space of the matrix $B(G)$.
\end{itemize}
\end{lem}
\noindent{\bf Proof.} According to the properties of orthogonal matrix and $B(G)B(G)^T=A(G)+D(G)$, then
\begin{itemize}
                        \item If $G$ is non-bipartite graph, then
\begin{eqnarray*}
  y_i^Ty_j=\left\{
             \begin{array}{ll}
               1, & \hbox{if $i=j$;} \\
               0, & \hbox{if $i\neq j$.}
             \end{array}
           \right.
\end{eqnarray*}

On the one hand,
\begin{eqnarray*}
  &&N(S_{k}(G))
\left(
                    \begin{array}{c}
                      v_i \\
                    \pm\frac{1}{\sqrt{k(1+\lambda_i)}}B^TD^{-\frac{1}{2}}v_i \\
                    \vdots \\
                    \pm\frac{1}{\sqrt{k(1+\lambda_i)}}B^TD^{-\frac{1}{2}}v_i
                    \end{array}
                  \right)
=\pm\sqrt{1+\lambda_i}\left(
                    \begin{array}{c}
                      v_i \\
                    \pm\frac{1}{\sqrt{k(1+\lambda_i)}}B^TD^{-\frac{1}{2}}v_i \\
                    \vdots \\
                    \pm\frac{1}{\sqrt{k(1+\lambda_i)}}B^TD^{-\frac{1}{2}}v_i
                    \end{array}
                  \right),\\
&&N(S_{k}(G))
\left(
                    \begin{array}{c}
                      0 \\
                    \sqrt{\frac{1}{k(k-1)(1+\lambda_i)}}B^TD^{-\frac{1}{2}}v_i \\
                    \sqrt{\frac{1}{k(k-1)(1+\lambda_i)}}B^TD^{-\frac{1}{2}}v_i \\
                    \sqrt{\frac{1}{k(k-1)(1+\lambda_i)}}B^TD^{-\frac{1}{2}}v_i\\
                    \sqrt{\frac{1}{k(k-1)(1+\lambda_i)}}B^TD^{-\frac{1}{2}}v_i \\
                    \vdots \\
                    -\sqrt{\frac{k-1}{k(1+\lambda_i)}}B^TD^{-\frac{1}{2}}v_i \\
                    \end{array}
                  \right)
=\left(
                    \begin{array}{c}
                      0 \\
                    0 \\
                    0 \\
                    0\\
                    0 \\
                    \vdots \\
                    0 \\
                    \end{array}
                  \right),
N(S_{k}(G))\left(
                    \begin{array}{c}
                      0 \\
                    y_z \\
                    0 \\
                    \vdots\\
                    0 \\
                    \end{array}
                  \right)
=\left(
                    \begin{array}{c}
                      0 \\
                    0 \\
                    0 \\
                    \vdots\\
                    0 \\
                    \end{array}
                  \right).
\end{eqnarray*}

On the other hand,
\begin{eqnarray*}
\left(
                    \begin{array}{c}
                      v_i \\
                    \pm\frac{1}{\sqrt{k(1+\lambda_i)}}B^TD^{-\frac{1}{2}}v_i \\
                    \vdots \\
                    \pm\frac{1}{\sqrt{k(1+\lambda_i)}}B^TD^{-\frac{1}{2}}v_i
                    \end{array}
                  \right)^T\left(
                    \begin{array}{c}
                      v_j \\
                    \pm\frac{1}{\sqrt{k(1+\lambda_j)}}B^TD^{-\frac{1}{2}}v_j \\
                    \vdots \\
                    \pm\frac{1}{\sqrt{k(1+\lambda_j)}}B^TD^{-\frac{1}{2}}v_j
                    \end{array}
                  \right)=\left\{
                            \begin{array}{ll}
                              2, & \hbox{if $i=j$ and the signs are the same;} \\
                              0, & \hbox{else.}
                            \end{array}
                          \right.
\end{eqnarray*}
\begin{eqnarray*}
\left(
                    \begin{array}{c}
                      0 \\
                    \sqrt{\frac{1}{k(k-1)(1+\lambda_i)}}B^TD^{-\frac{1}{2}}v_i \\
                    \sqrt{\frac{1}{k(k-1)(1+\lambda_i)}}B^TD^{-\frac{1}{2}}v_i \\
                    \sqrt{\frac{1}{k(k-1)(1+\lambda_i)}}B^TD^{-\frac{1}{2}}v_i\\
                    \sqrt{\frac{1}{k(k-1)(1+\lambda_i)}}B^TD^{-\frac{1}{2}}v_i \\
                    \vdots \\
                    -\sqrt{\frac{k-1}{k(1+\lambda_i)}}B^TD^{-\frac{1}{2}}v_i \\
                    \end{array}
                  \right)^T\left(
                    \begin{array}{c}
                      0 \\
                    \sqrt{\frac{1}{k(k-1)(1+\lambda_j)}}B^TD^{-\frac{1}{2}}v_j \\
                    \sqrt{\frac{1}{k(k-1)(1+\lambda_j)}}B^TD^{-\frac{1}{2}}v_j \\
                    \sqrt{\frac{1}{k(k-1)(1+\lambda_j)}}B^TD^{-\frac{1}{2}}v_j\\
                    \sqrt{\frac{1}{k(k-1)(1+\lambda_j)}}B^TD^{-\frac{1}{2}}v_j \\
                    \vdots \\
                    -\sqrt{\frac{k-1}{k(1+\lambda_j)}}B^TD^{-\frac{1}{2}}v_j \\
                    \end{array}
                  \right)=\left\{
                            \begin{array}{ll}
                              1, & \hbox{if $i=j$;} \\
                              0, & \hbox{else.}
                            \end{array}
                          \right.
\end{eqnarray*}
\begin{eqnarray*}
\left(
                    \begin{array}{c}
                      0 \\
                    y_i \\
                    0 \\
                    \vdots\\
                    0 \\
                    \end{array}
                  \right)^T\left(
                    \begin{array}{c}
                      0 \\
                    y_j \\
                    0 \\
                    \vdots\\
                    0 \\
                    \end{array}
                  \right)=\left\{
                            \begin{array}{ll}
                              1, & \hbox{if $i=j$;} \\
                              0, & \hbox{else.}
                            \end{array}
                          \right.
\end{eqnarray*}

Moreover, one obtains that
\begin{eqnarray*}
\left(
                    \begin{array}{c}
                      v_i \\
                    \pm\frac{1}{\sqrt{k(1+\lambda_i)}}B^TD^{-\frac{1}{2}}v_i \\
                    \vdots \\
                    \pm\frac{1}{\sqrt{k(1+\lambda_i)}}B^TD^{-\frac{1}{2}}v_i
                    \end{array}
                  \right)^T\left(
                    \begin{array}{c}
                      0 \\
                    \sqrt{\frac{1}{k(k-1)(1+\lambda_i)}}B^TD^{-\frac{1}{2}}v_i \\
                    \sqrt{\frac{1}{k(k-1)(1+\lambda_i)}}B^TD^{-\frac{1}{2}}v_i \\
                    \sqrt{\frac{1}{k(k-1)(1+\lambda_i)}}B^TD^{-\frac{1}{2}}v_i\\
                    \sqrt{\frac{1}{k(k-1)(1+\lambda_i)}}B^TD^{-\frac{1}{2}}v_i \\
                    \vdots \\
                    -\sqrt{\frac{k-1}{k(1+\lambda_i)}}B^TD^{-\frac{1}{2}}v_i \\
                    \end{array}
                  \right)=0,
\end{eqnarray*}
\begin{eqnarray*}
\left(
                    \begin{array}{c}
                      v_i \\
                    \pm\frac{1}{\sqrt{k(1+\lambda_i)}}B^TD^{-\frac{1}{2}}v_i \\
                    \vdots \\
                    \pm\frac{1}{\sqrt{k(1+\lambda_i)}}B^TD^{-\frac{1}{2}}v_i
                    \end{array}
                  \right)^T\left(
                    \begin{array}{c}
                      0 \\
                    y_z \\
                    0 \\
                    \vdots\\
                    0 \\
                    \end{array}
                  \right)=0,~
\left(
                    \begin{array}{c}
                      0 \\
                    \sqrt{\frac{1}{k(k-1)(1+\lambda_i)}}B^TD^{-\frac{1}{2}}v_i \\
                    \sqrt{\frac{1}{k(k-1)(1+\lambda_i)}}B^TD^{-\frac{1}{2}}v_i \\
                    \sqrt{\frac{1}{k(k-1)(1+\lambda_i)}}B^TD^{-\frac{1}{2}}v_i\\
                    \sqrt{\frac{1}{k(k-1)(1+\lambda_i)}}B^TD^{-\frac{1}{2}}v_i \\
                    \vdots \\
                    -\sqrt{\frac{k-1}{k(1+\lambda_i)}}B^TD^{-\frac{1}{2}}v_i \\
                    \end{array}
                  \right)^T\left(
                    \begin{array}{c}
                      0 \\
                    y_z \\
                    0 \\
                    \vdots\\
                    0 \\
                    \end{array}
                  \right)=0,
\end{eqnarray*}
Hence, those orthonormal eigenvectors of the corresponding eigenvalues for the graph $S_{k}(G)$ as shown in theorem.
 \item If $G$ is non-bipartite graph, the proofs of this case are similar to the above procedure. Hence, we omit here.
                      \end{itemize}
The result as desired. \hfill\rule{1ex}{1ex}

The Lemma 4.5 has given the orthonormal eigenvectors of the corresponding eigenvalues in terms of Theorem 1.1 of $N(S_{k}(G))$. At this point, we will put the orthonormal eigenvectors of $N(S_{k}(G))$ in another way.
\begin{enumerate}
\item The orthonormal eigenvectors of the corresponding eigenvalues $\pm\sqrt{\frac{1+\lambda_1}{2}}=\pm 1$ are as follows.
\begin{eqnarray*}
\bigg(\sqrt{\frac{kd_1}{4km}},  \sqrt{\frac{kd_2}{4km}},  \ldots,  \sqrt{\frac{kd_n}{4km}},  \sqrt{\frac{2}{4km}},
     \sqrt{\frac{2}{4km}},  \ldots,  \sqrt{\frac{2}{4km}} \bigg)^T,\\
\bigg(\sqrt{\frac{kd_1}{4km}},  \sqrt{\frac{kd_2}{4km}},  \ldots,  \sqrt{\frac{kd_n}{4km}},  -\sqrt{\frac{2}{4km}},
     -\sqrt{\frac{2}{4km}},  \ldots,  -\sqrt{\frac{2}{4km}} \bigg)^T,
\end{eqnarray*}
respectively.
\item The eigenvectors $u_j$ of the corresponding eigenvalue $\pm\sqrt{\frac{1+\lambda_a}{2}},~a=2,3,\ldots,n$ are as follows.
\begin{eqnarray*}
u_j=\left\{
  \begin{array}{ll}
    \frac{\sqrt{2}}{2}v_{kj}, & \hbox{if $j\in V(G)$;} \\
    \pm\sqrt{\frac{1}{2k(1+\lambda_a)}}\big(\frac{v_{as}}{\sqrt{d_s}}+\frac{v_{at}}{\sqrt{d_t}}\big), & \hbox{if $j\in N$ with $N_{S_{k}(G)}(j)=\{s,t\}$.}
  \end{array}
\right.
\end{eqnarray*}
Based on the properties of the orthogonal matrix, yields
\begin{eqnarray*}
\sum_{l=1}^{m-n}y_{lj}^2&=&1-\frac{1}{km}-\sum_{a=2}^{n}\frac{1}{k(1+\lambda_a)}\bigg(\frac{v_{as}}{\sqrt{d_s}}+\frac{v_{at}}{\sqrt{d_t}}\bigg)^2
-\sum_{a=1}^{n}\frac{k-1}{k(1+\lambda_a)}\bigg(\frac{v_{as}}{\sqrt{d_s}}+\frac{v_{at}}{\sqrt{d_t}}\bigg)^2\\
&=&1-\frac{1}{m}-\sum_{a=2}^{n}\frac{1}{1+\lambda_a}\bigg(\frac{v_{as}}{\sqrt{d_s}}+\frac{v_{at}}{\sqrt{d_t}}\bigg)^2,
\end{eqnarray*}
if $G$ is non-bipartite graph.
\item The eigenvectors $u_j$ of the corresponding eigenvalue $\pm\sqrt{\frac{1+\lambda_a}{2}},~a=2,3,\ldots,n-1$ are as follows.
\begin{eqnarray*}
u_j=\left\{
  \begin{array}{ll}
    \frac{\sqrt{2}}{2}v_{kj}, & \hbox{if $j\in V(G)$;} \\
    \pm\sqrt{\frac{1}{2k(1+\lambda_a)}}\big(\frac{v_{as}}{\sqrt{d_s}}+\frac{v_{at}}{\sqrt{d_t}}\big), & \hbox{if $j\in N$ with $N_{S_{k}(G)}(j)=\{s,t\}$.}
  \end{array}
\right.
\end{eqnarray*}
Based on the properties of the orthogonal matrix, yields
\begin{eqnarray*}
\sum_{l=1}^{m-n+1}y_{lj}^2&=&1-\frac{1}{km}-\sum_{a=2}^{n-1}\frac{1}{k(1+\lambda_a)}\bigg(\frac{v_{as}}{\sqrt{d_s}}+\frac{v_{at}}{\sqrt{d_t}}\bigg)^2
-\sum_{a=1}^{n-1}\frac{k-1}{k(1+\lambda_a)}\bigg(\frac{v_{as}}{\sqrt{d_s}}+\frac{v_{at}}{\sqrt{d_t}}\bigg)^2\\
&=&1-\frac{1}{m}-\sum_{a=2}^{n-1}\frac{1}{1+\lambda_a}\bigg(\frac{v_{as}}{\sqrt{d_s}}+\frac{v_{at}}{\sqrt{d_t}}\bigg)^2,
\end{eqnarray*}
if $G$ is bipartite graph.
\end{enumerate}

According to the structures of the graph $S_{k}(G)$, the selections of vertices $i$ and $j$ can be divided into three cases while determine the hitting time between vertices $i$ and $j$.

\noindent{\bf Case 1.} $i,~j\in V(G)$ and $G$ is non-bipartite graph, then
\begin{eqnarray*}
E_iT_j(S_{k}(G))&=&4km\sum_{a=2}^{n}\bigg(\frac{1}{1-\sqrt{\frac{1+\lambda_a}{2}}}+\frac{1}{1+\sqrt{\frac{1+\lambda_a}{2}}}\bigg)
\bigg(\frac{v_{aj}^2}{2kd_j}+\frac{v_{aj}v_{ai}}{2k\sqrt{d_id_j}}\bigg)\\
&=&4km\sum_{a=2}^{n}\frac{4}{1-\lambda_a}
\bigg(\frac{v_{aj}^2}{2kd_j}+\frac{v_{aj}v_{ai}}{2k\sqrt{d_id_j}}\bigg)\\
&=&8m\sum_{a=2}^{n}\frac{1}{1-\lambda_a}\bigg(\frac{v_{aj}^2}{d_j}+\frac{v_{aj}v_{ai}}{\sqrt{d_id_j}}\bigg)\\
&=&4E_iT_j(G),
\end{eqnarray*}
while $G$ is bipartite graph, then
\begin{eqnarray*}
E_iT_j(S_{k}(G))&=&4km\Bigg[\sum_{a=2}^{n-1}\bigg(\frac{1}{1-\sqrt{\frac{1+\lambda_a}{2}}}+\frac{1}{1+\sqrt{\frac{1+\lambda_a}{2}}}\bigg)
\bigg(\frac{v_{aj}^2}{2kd_j}+\frac{v_{aj}v_{ai}}{2k\sqrt{d_id_j}}\bigg)+\frac{v_{nj}^2}{kd_j}+\frac{v_{nj}v_{ni}}{k\sqrt{d_id_j}}\Bigg]\\
&=&4km\Bigg[\sum_{a=2}^{n-1}\frac{2}{k}\cdot\frac{1}{1-\lambda_a}
\bigg(\frac{v_{aj}^2}{d_j}+\frac{v_{aj}v_{ai}}{\sqrt{d_id_j}}\bigg)+\frac{v_{nj}^2}{kd_j}+\frac{v_{nj}v_{ni}}{k\sqrt{d_id_j}}\Bigg]\\
&=&8m\sum_{a=2}^{n}\frac{1}{1-\lambda_a}\bigg(\frac{v_{aj}^2}{d_j}+\frac{v_{aj}v_{ai}}{\sqrt{d_id_j}}\bigg)\\
&=&4E_iT_j(G).
\end{eqnarray*}
Hence, $E_iT_j(S_{k}(G))=4E_iT_j(G)$ holds whatever $G$ is bipartite graph or not.

\noindent{\bf Case 2.} $i\in N,~j\in V(G)$ and $N_{S_k(G)}(i)=\{s,t\}$, then
\begin{eqnarray*}
E_iT_j(S_{k}(G))=1+\frac{1}{2}\big(E_sT_j(S_k(G))+E_tT_j(S_k(G))\big)=1+2E_sT_j(G)+2E_tT_j(G).
\end{eqnarray*}

For $E_jT_i(S_{k}(G))$, one obtains the following equation based on Lemma 2.5, Enumerates 1 and 2.
\begin{eqnarray*}
E_jT_i(S_{k}(G))&=&4km\Bigg[\frac{1}{4km}+\sum_{a=2}^{n}\bigg(\frac{1}{2-2\sqrt{\frac{1+\lambda_a}{2}}}+\frac{1}{2+2
\sqrt{\frac{1+\lambda_a}{2}}}\bigg)
\frac{1}{2k(1+\lambda_a)}\bigg(\frac{v_{as}}{\sqrt{d_s}}+\frac{v_{at}}{\sqrt{d_t}}\bigg)^2\\
&&-\sum_{a=2}^{n}\bigg(\frac{1}{2k-2k\sqrt{\frac{1+\lambda_a}{2}}}-\frac{1}{2k+2k\sqrt{\frac{1+\lambda_a}{2}}}\bigg)
\frac{v_{aj}}{\sqrt{2(1+\lambda_a)d_j}}\bigg(\frac{v_{as}}{\sqrt{d_s}}+\frac{v_{at}}{\sqrt{d_t}}\bigg)\\
&&+\sum_{l=1}^{m-n}\frac{y_{li}^2}{2}+\sum_{a=1}^{n}\frac{k-1}{2k(1+\lambda_a)}\bigg(\frac{v_{as}}{\sqrt{d_s}}+
\frac{v_{at}}{\sqrt{d_t}}\bigg)^2\Bigg]\\
&=&4km\Bigg[\sum_{a=2}^{n}\frac{1}{k(1+\lambda_a)(1-\lambda_a)}
\bigg(\frac{v_{as}}{\sqrt{d_s}}+\frac{v_{at}}{\sqrt{d_t}}\bigg)^2-\sum_{a=2}^{n}\frac{1}{k(1-\lambda_a)}
\bigg(\frac{v_{as}v_{aj}}{\sqrt{d_sd_j}}+\frac{v_{at}v_{aj}}{\sqrt{d_td_j}}\bigg)\\
&&+\frac{1}{4km}+1-\frac{1}{m}-\sum_{a=2}^{n}\frac{1}{1+\lambda_a}\bigg(\frac{v_{as}}{\sqrt{d_s}}+\frac{v_{at}}{\sqrt{d_t}}\bigg)^2
+\sum_{a=1}^{n}\frac{k-1}{2k(1+\lambda_a)}\bigg(\frac{v_{as}}{\sqrt{d_s}}+\frac{v_{at}}{\sqrt{d_t}}\bigg)^2\Bigg]\\
&=&2km-1+4m\sum_{a=2}^{n}\frac{1}{1-\lambda_a}\bigg[\frac{1}{2}\bigg(\frac{v_{as}}{\sqrt{d_s}}+\frac{v_{at}}{\sqrt{d_t}}\bigg)^2
-\frac{v_{as}v_{aj}}{\sqrt{d_sd_j}}-\frac{v_{at}v_{aj}}{\sqrt{d_td_j}}\bigg]\\
&=&2km-1+4m\sum_{a=2}^{n}\frac{1}{1-\lambda_a}\bigg(\frac{v_{as}^2}{d_s}+\frac{v_{at}^2}{d_t}
-\frac{v_{as}v_{aj}}{\sqrt{d_sd_j}}-\frac{v_{at}v_{aj}}{\sqrt{d_td_j}}-\frac{v_{as}^2}{2d_s}-\frac{v_{at}^2}{2d_t}+
\frac{v_{as}v_{at}}{\sqrt{d_sd_t}}\bigg)\\
&=&2km-1+2\big[E_jT_s(G)+E_jT_t(G)\big]-\big[E_tT_s(G)+E_sT_t(G)\big],
\end{eqnarray*}
if $G$ is non-bipartite graph.

In the same way, one arrives at
\begin{eqnarray*}
E_jT_i(S_{k}(G))&=&4km\Bigg[\frac{1}{4km}+\sum_{a=2}^{n}\bigg(\frac{1}{2-2\sqrt{\frac{1+\lambda_a}{2}}}+\frac{1}{2+2
\sqrt{\frac{1+\lambda_a}{2}}}\bigg)
\frac{1}{2k(1+\lambda_a)}\bigg(\frac{v_{as}}{\sqrt{d_s}}+\frac{v_{at}}{\sqrt{d_t}}\bigg)^2\\
&&-\sum_{a=2}^{n-1}\bigg(\frac{1}{2k-2k\sqrt{\frac{1+\lambda_a}{2}}}-\frac{1}{2k+2k\sqrt{\frac{1+\lambda_a}{2}}}\bigg)
\frac{v_{aj}}{\sqrt{2(1+\lambda_a)d_j}}\bigg(\frac{v_{as}}{\sqrt{d_s}}+\frac{v_{at}}{\sqrt{d_t}}\bigg)\\
&&+\sum_{l=1}^{m-n+1}\frac{y_{li}^2}{2}+\sum_{a=1}^{n}\frac{k-1}{2k(1+\lambda_a)}\bigg(\frac{v_{as}}{\sqrt{d_s}}+
\frac{v_{at}}{\sqrt{d_t}}\bigg)^2\Bigg]\\
%&=&4km\Bigg[\sum_{a=2}^{n}\frac{1}{k(1+\lambda_a)(1-\lambda_a)}
%\bigg(\frac{v_{as}}{\sqrt{d_s}}+\frac{v_{at}}{\sqrt{d_t}}\bigg)^2-\sum_{a=2}^{n}\frac{1}{k(1-\lambda_a)}
%\bigg(\frac{v_{as}v_{aj}}{\sqrt{d_sd_j}}+\frac{v_{at}v_{aj}}{\sqrt{d_td_j}}\bigg)\\
%&&+\frac{1}{4km}+1-\frac{1}{m}-\sum_{a=2}^{n}\frac{1}{1+\lambda_a}\bigg(\frac{v_{as}}{\sqrt{d_s}}+\frac{v_{at}}{\sqrt{d_t}}\bigg)^2
%+\sum_{a=1}^{n}\frac{k-1}{2k(1+\lambda_a)}\bigg(\frac{v_{as}}{\sqrt{d_s}}+\frac{v_{at}}{\sqrt{d_t}}\bigg)^2\Bigg]\\
%&=&2km-1+4m\sum_{a=2}^{n}\frac{1}{1-\lambda_a}\bigg[\frac{1}{2}\bigg(\frac{v_{as}}{\sqrt{d_s}}+\frac{v_{at}}{\sqrt{d_t}}\bigg)^2
%-\frac{v_{as}v_{aj}}{\sqrt{d_sd_j}}-\frac{v_{at}v_{aj}}{\sqrt{d_td_j}}\bigg]\\
&=&2km-1+4m\sum_{a=2}^{n-1}\frac{1}{1-\lambda_a}\bigg(\frac{v_{as}^2}{d_s}+\frac{v_{at}^2}{d_t}
-\frac{v_{as}v_{aj}}{\sqrt{d_sd_j}}-\frac{v_{at}v_{aj}}{\sqrt{d_td_j}}-\frac{v_{as}^2}{2d_s}-\frac{v_{at}^2}{2d_t}+
\frac{v_{as}v_{at}}{\sqrt{d_sd_t}}\bigg)\\
&=&2km-1+2\big[E_jT_s(G)+E_jT_t(G)\big]-\big[E_tT_s(G)+E_sT_t(G)\big],
\end{eqnarray*}
if $G$ is bipartite graph.

\noindent{\bf Case 3.} $i,~j\in N$, $N_{S_k(G)}(i)=\{s,t\}$ and $N_{S_k(G)}(j)=\{p,q\}$. According to Case 2, then
\begin{eqnarray*}
E_iT_j(S_{k}(G))&=&1+\frac{1}{2}\big[E_sT_j(S_k(G))+E_tT_j(S_k(G))\big]\\
&=&1+\frac{1}{2}\big[2km-1+2\big(E_sT_p(G)+E_sT_q(G)\big)-\big(E_qT_p(G)+E_pT_q(G)\big)\big]\\
&&+\frac{1}{2}\big[2km-1+2\big(E_tT_p(G)+E_tT_q(G)\big)-\big(E_qT_p(G)+E_pT_q(G)\big)\big]\\
&=&2km+E_sT_p(G)+E_sT_q(G)+E_tT_p(G)+E_tT_q(G)-E_qT_p(G)-E_pT_q(G).
\end{eqnarray*}

In the same way, one obtains
\begin{eqnarray*}
E_jT_i(S_{k}(G))=2km+E_pT_s(G)+E_qT_s(G)+E_pT_t(G)+E_qT_t(G)-E_sT_p(G)-E_pT_s(G).
\end{eqnarray*}

Combining Cases 1, 2 and 3, we get the following theorem.
\begin{thm} Assume that $G$ is a $(n,m)$-graph, then the expected hitting time between vertices $i$ and $j$ in $S_k(G)$ are as follows.
\begin{itemize}
  \item $i,~j\in V(G)$, then~$$E_iT_j(S_{k}(G))=4E_iT_j(G).$$
  \item $i\in N,~j\in V(G)$ and $N_{S_k(G)}(i)=\{s,t\}$, then
\begin{eqnarray*}
E_iT_j(S_{k}(G))=1+2E_sT_j(G)+2E_tT_j(G),
\end{eqnarray*}
\begin{eqnarray*}
E_jT_i(S_{k}(G))=2km-1+2\big[E_jT_s(G)+E_jT_t(G)\big]-\big[E_tT_s(G)+E_sT_t(G)\big].
\end{eqnarray*}
  \item $i,~j\in N$, $N_{S_k(G)}(i)=\{s,t\}$ and $N_{S_k(G)}(j)=\{p,q\}$, then
\begin{eqnarray*}
E_iT_j(S_{k}(G))=2km+E_sT_p(G)+E_sT_q(G)+E_tT_p(G)+E_tT_q(G)-E_qT_p(G)-E_pT_q(G),
\end{eqnarray*}
\begin{eqnarray*}
E_jT_i(S_{k}(G))=2km+E_pT_s(G)+E_qT_s(G)+E_pT_t(G)+E_qT_t(G)-E_sT_p(G)-E_pT_s(G).
\end{eqnarray*}
\end{itemize}
\end{thm}

Based on Lemma 2.6 and above theorem, leads
\begin{cor} Assume that $G$ is a $(n,m)$-graph, then any two-points resistance distance between vertices $i$ and $j$ in $S_k(G)$ are as below.
\begin{itemize}
  \item $i,~j\in V(G)$, then~$$\Omega_{ij}(S_{k}(G))=\frac{2}{k}\Omega_{ij}(G).$$
  \item $i\in N,~j\in V(G)$ and $N_{S_k(G)}(i)=\{s,t\}$, then
\begin{eqnarray*}
\Omega_{ij}(S_{k}(G))=\frac{k+2\Omega_{sj}(G)+2\Omega_{tj}(G)-\Omega_{st}(G)}{2k}.
\end{eqnarray*}
  \item $i,~j\in N$, $N_{S_k(G)}(i)=\{s,t\}$ and $N_{S_k(G)}(j)=\{p,q\}$, then
\begin{eqnarray*}
\Omega_{ij}(S_{k}(G))=\frac{2k+\Omega_{sp}(G)+\Omega_{sq}(G)+\Omega_{tp}(G)+\Omega_{tq}(G)-\Omega_{pq}(G)-\Omega_{st}(G)}{2k}.
\end{eqnarray*}
\end{itemize}
\end{cor}
\noindent{\bf Remark 1.} While $k=2$, we can get the Theorem 4.2 of\cite{Guo}. Indeed, one can compute the multiplicative degree-Kirchhoff index of $S_k(G)$ through any two-points resistance distance between vertices $i$ and $j$ in $S_k(G)$. But in this paper, we prefer to use the normalized Laplacian spectra of $S_k(G)$ to calculate the multiplicative degree-Kirchhoff index of $S_k(G)$(see Eq.(4.20)). This seems much easier.

According to the relation between the expected hitting time and expected commute time, the following corollary can be immediately obtained.
\begin{cor} Assume that $G$ is a $(n,m)$-graph, then the expected commute time between vertices $i$ and $j$ in $S_k(G)$ are as below.
\begin{itemize}
  \item $i,~j\in V(G)$, then~$$C_{ij}(S_{k}(G))=\frac{4m}{k}\Omega_{ij}(G).$$
  \item $i\in N,~j\in V(G)$ and $N_{S_k(G)}(i)=\{s,t\}$, then
\begin{eqnarray*}
C_{ij}(S_{k}(G))=m\cdot\frac{k+2C_{sj}(G)+2C_{tj}(G)-C_{st}(G)}{k}.
\end{eqnarray*}
  \item $i,~j\in N$, $N_{S_k(G)}(i)=\{s,t\}$ and $N_{S_k(G)}(j)=\{p,q\}$, then
\begin{eqnarray*}
C_{ij}(S_{k}(G))=m\cdot\frac{2k+C_{sp}(G)+C_{sq}(G)+C_{tp}(G)+C_{tq}(G)-C_{pq}(G)-C_{st}(G)}{k}.
\end{eqnarray*}
\end{itemize}
\end{cor}

It remains to now look at the relation between the normalized Laplacians of $G$ and $S_k(G)$. First goal is to determine the eigenvalues of $\mathscr{L}(S_k(G))$ with $r$ iterations, i.e., the eigenvalues of $\mathscr{L}(S_k^r(G))$. Second, we provide in this part
the explicit formulas for $Kf^*(S_{k}(G))$($Kf^*(S_k^r(G))$ resp.), $Ke(S_{k}(G))$($Ke(S_k^r(G))$ resp.) and $\tau(S_{k}(G))$($\tau(S_k^r(G))$ resp.).

Let $S_k^0(G)=G$, $S_k^1(G)=S_k(S_k^0(G))$, \ldots, $S_k^{r}(G)=S_k(S_k^{r-1}(G))$. Then denote $|E(S_k^{r}(G))|$ and $|V(S_k^{r}(G))|$ the edge set and vertex set of $S_k^r(G)$, $|E_r|$ and $|V_r|$ for simplify. According to the construction of $S_k^r(G)$, one has
\begin{eqnarray*}
|E_r|=2k|E_{r-1}|,~|V_r|=|V_{r-1}|+k|E_{r-1}|.
\end{eqnarray*}
It is easy to obtain
\begin{eqnarray*}
|E_r|=m(2k)^r,~|V_r|=n+km\frac{(2k)^r-1}{2k-1}.
\end{eqnarray*}

Assume that $U$ is a finite multiset of real numbers. To obtain the normalized Laplacians of $S_k^r(G)$, we also have to define two new multisets in the line with Eq.(3.14) as follows.
\begin{eqnarray*}
f_1(U)=\bigcup_{x\in U}\{f_1(x)\},~~~f_2(U)=\bigcup_{x\in U}\{f_2(x)\}.
\end{eqnarray*}

With those notation in hand, the normalized Laplacians of $S_k^r(G)$ can be determined immediately by Theorem 1.1 and the structure of $S_k^r(G)$.
\begin{thm} Assume that $S_k^r(G)$ is the r-th iterations of the graph $S_k(G)$, then the normalized Laplacian spectra of $S_k^r(G)$ are as follows.
\begin{eqnarray*}
\begin{cases}
f_1(\Gamma(G)\backslash \{0,2\})\bigcup f_2(\Gamma(G)\backslash \{0,2\})\bigcup \{0,2\}\bigcup\{\underbrace{1,1, \ldots, 1}_{km-n+2}\},~~~r=1~and~G~is~bipartite;\\
f_1(\Gamma(G))\bigcup f_2(\Gamma(G))\bigcup \{\underbrace{1,1, \ldots, 1}_{km-n}\},~~~r=1~and~G~is~nonbipartite;\\
f_1(\Gamma(S_k^{r-1}(G))\backslash \{0,2\})\bigcup f_2(\Gamma(S_k^{r-1}(G))\backslash \{0,2\})\bigcup \{0,2\}\bigcup\{1,1, \ldots, 1\},~~~r>1,\\
\end{cases}
\end{eqnarray*}
where $\Gamma(S_k^{r-1}(G))$ is the normalized Laplacian spectrum of $S_k^{r-1}(G)$, the multiplicities are respectively $2^{r-1}(n-1)+\sum_{j=2}^r2^{j-2}(k|E_{r-j}|-|V_{r-j}|+2)$ and $k|E_{r-1}|-|V_{r-1}|+2$ of $f_1(\Gamma(S_k^{r-1}(G))\backslash \{0,2\})$ and $\{1,1, \ldots, 1\}$ for $r>1$.
\end{thm}

According to the normalized Laplacians of $S_{k}^r(G)$, we derive the explicit formulas for $Kf^*(S_k^r(G))$, $Ke(S_k^r(G))$ and $\tau(S_k^r(G))$.
\begin{thm}
Assume that $G$ is a $(n,m)$-graph. One has
\begin{eqnarray*}
Kf^*(S_k^r(G))&=&(8k)^rKf^*(G)+\frac{(2k)^r(4^r-1)}{3}(m-2mn)+\frac{k(4k)^r(k^r-2^r)}{k-2}m^2\\
&&-\frac{k(2k)^r[4^r-2k(4^r-1)+3(2k)^r-4]}{3(k-2)(2k-1)}m^2,
\end{eqnarray*}
where $k> 2,~r\geq 1$.
\end{thm}
\noindent{\bf Proof.} According to Theorem 1.1 and $\lambda_i\neq 2$, one has $f_1(\lambda_i)$ and $f_2(\lambda_i)$ are two roots of $-2\sigma^2+4\sigma-\lambda_i=0$, then
\begin{eqnarray}
\frac{1}{f_1(\lambda_i)}+\frac{1}{f_2(\lambda_i)}=\frac{4}{\lambda_i},~f_1(\lambda_i)f_2(\lambda_i)=\frac{\lambda_i}{2}.
\end{eqnarray}
In addition, $\frac{1}{1}+\frac{1}{1}=2,~1\times 1=\frac{2}{2}=1$, thus $\lambda_i= 2$ also satisfy Eq.(4.19). Based on Lemma (2.1), one obtains
\begin{eqnarray}
Kf^*(S_k^1(G))&=&4km\bigg(\sum_{i=2}^n\bigg(\frac{1}{f_1(\lambda_i)}+\frac{1}{f_2(\lambda_i)}\bigg)+\frac{1}{2}+km-n\bigg)\nonumber\\
&=&4km\bigg(\sum_{i=2}^n\frac{4}{\lambda_i}+\frac{1}{2}+km-n\bigg)\nonumber\\
&=&8k\cdot Kf^*(G)+2km(1+2km-2n).
\end{eqnarray}

By Eq.(4.20) and the construction of $S_{k}^r(G)$, we get
\begin{eqnarray*}
Kf^*(S_k^r(G))&=&8k\cdot Kf^*(S_k^{r-1}(G))+2k|E_{r-1}|(1+2k|E_{r-1}|-2|V_{r-1}|)\\
&=&(8k)^rKf^*(G)+2k\sum_{i=0}^{r-1}(8k)^{r-1-i}|E_i|(1+2k|E_i|-2|V_i|)\\
&=&(8k)^rKf^*(G)+\frac{(2k)^r(4^r-1)}{3}(m-2mn)+\frac{k(4k)^r(k^r-2^r)}{k-2}m^2\\
&&-\frac{k(2k)^r[4^r-2k(4^r-1)+3(2k)^r-4]}{3(k-2)(2k-1)}m^2.
\end{eqnarray*}
This completes the proof.\hfill\rule{1ex}{1ex}

\noindent{\bf Remark 2.} When $k=1$, it is
\begin{eqnarray*}
Kf^*(S(G))=8\cdot Kf^*(G)+2m(1+2m-2n).
\end{eqnarray*}
This coincides the Theorem 4.1 of \cite{.P}.

\noindent{\bf Remark 3.} When $k=2$, it satisfies
\begin{eqnarray*}
Kf^*(S_2^1(G))=16\cdot Kf^*(G)+4m(1+4m-2n).
\end{eqnarray*}
This equals the Theorem 4.1 of\cite{Guo}.

For $S_2^r(G)$, the multiplicative degree-Kirchhoff index is as below.
\begin{eqnarray*}
Kf^*(S_2^r(G))&=&16\cdot Kf^*(S_2^{r-1}(G))+4|E_{r-1}|(1+4|E_{r-1}|-2|V_{r-1}|)\\
&=&(16)^rKf^*(G)+4\sum_{i=0}^{r-1}(16)^{r-1-i}|E_i|(1+4|E_i|-2|V_i|)\\
&=&(16)^rKf^*(G)+\frac{4^r(4^r-1)}{3}m-\frac{2\cdot 4^r(4^r-1)}{3}mn\\
&&+\frac{2\cdot4^r[2(4^r-1)+3r\cdot4^r]}{9}m^2.
\end{eqnarray*}

\begin{thm}
Assume that $G$ is a $(n,m)$-graph. One has
\begin{eqnarray*}
Ke(S_k^r(G))=4^rKe(G)+\frac{4^r-1}{6}(1-2n)+\frac{2^{r-1}k(k^r-2^r)}{k-2}m-\frac{k[4^r-2k(4^r-1)+3(2k)^r-4]}{6(k-2)(2k-1)}m,
\end{eqnarray*}
where $k> 2,~r\geq 1$.
\end{thm}
\noindent{\bf Proof.} By Lemma (2.2) and the relation between $Ke(G)$ and $Kf^*(G)$, one gets
\begin{eqnarray*}
Ke(S_k^1(G))=\frac{1}{4km}Kf^*(S_k^1(G))=4Ke(G)+\frac{1}{2}(1+2km-2n).
\end{eqnarray*}
Further,
\begin{eqnarray*}
Ke(S_k^r(G))&=&\frac{1}{2m(2k)^r}Kf^*(S_k^r(G))\\
&=&\frac{(8k)^r}{2m(2k)^r}Kf^*(G)+\frac{(2k)^r(4^r-1)}{6m(2k)^r}(m-2mn)+\frac{k(4k)^r(k^r-2^r)}{2m(2k)^r(k-2)}m^2\\
&&-\frac{k(2k)^r[4^r-2k(4^r-1)+3(2k)^r-4]}{6m(2k)^r(k-2)(2k-1)}m^2\\
&=&4^rKe(G)+\frac{4^r-1}{6}(1-2n)+\frac{2^{r-1}k(k^r-2^r)}{k-2}m-\frac{k[4^r-2k(4^r-1)+3(2k)^r-4]}{6(k-2)(2k-1)}m.
\end{eqnarray*}
The proof has completed.\hfill\rule{1ex}{1ex}

\noindent{\bf Remark 4.} When $k=2$, this yields
\begin{eqnarray*}
Ke(S_2^r(G))&=&\frac{1}{2m\cdot4^r}Kf^*(S_2^r(G))\\
&=&4^rKe(G)+\frac{4^r-1}{6}(1-2n)+\frac{2(4^r-1)+3r\cdot4^r}{9}m.
\end{eqnarray*}
While $r=1$, we can directly obtain the theorem 4.1 in\cite{Guo}.

In the latest part of this section, the formula for $\tau(S_{k}^r(G))$ is determined.
\begin{thm}
Assume that $G$ is a $(n,m)$-graph. One has
\begin{eqnarray*}
\tau(S_{k}^r(G))=2^{km\cdot\frac{1-(2k)^r}{1-2k}-r}\cdot k^{nr+km\cdot\frac{(2k)^r-2kr+r-1}{(1-2k)^2}-r}\tau(G),
\end{eqnarray*}
where $k\geq 1,~r\geq 1$.
\end{thm}
\noindent{\bf Proof.} Based on Eq.(3.2), it is routine to check
\begin{eqnarray}
\prod_{i=1}^{n+km}d_i(S_k(G))=2^{km}\cdot\prod_{i=1}^{n}kd_i(G)=2^{km}\cdot k^n\cdot\prod_{i=1}^{n}d_i(G).
\end{eqnarray}
By Eq.(4.21) and Lemma (2.4), we obtain
\begin{eqnarray*}
\tau(S_{k}^1(G))&=&\frac{2\prod_{i=1}^{n+km}d_i(S_k(G))\prod^n_{i=2}f_1(\lambda_i)f_2(\lambda_i)}{4km}\\
&=&\frac{2^{km}\cdot k^n\cdot\prod_{i=1}^{n}d_i(G)\prod^n_{i=2}\lambda_i}{4km}\\
&=&2^{km-1}\cdot k^{n-1}\cdot\tau(G).
\end{eqnarray*}

Hence, one arrives at
\begin{eqnarray*}
\tau(S_{k}^r(G))&=&2^{k|E_{r-1}|-1}\cdot k^{|V_{r-1}|-1}\cdot\tau(S_k^{r-1}(G))\\
&=&2^{\sum_{i=0}^{r-1}(k|E_i|-1)}\cdot k^{\sum_{i=0}^{r-1}(|V_i|-1)}\tau(G)\\
&=&2^{km\cdot\frac{1-(2k)^r}{1-2k}-r}\cdot k^{nr+km\cdot\frac{(2k)^r-2kr+r-1}{(1-2k)^2}-r}\tau(G).
\end{eqnarray*}
The result as desired. \hfill\rule{1ex}{1ex}

\section{Applications of Theorem 1.2}
\ \ \ \
Let $S_{2k}^0(G)=G$, $S_{2k}^1(G)=S_{2k}(S_{2k}^0(G))$, \ldots, $S_{2k}^{r}(G)=S_{2k}(S_{2k}^{r-1}(G))$. Then denote $|E(S_{2k}^{r}(G))|$ and $|V(S_{2k}^{r}(G))|$ the edge set and vertex set of $S_{2k}^r(G)$, $|E'_r|$ and $|V'_r|$ for simplify. According to the construction of $S_{2k}^r(G)$, one has
\begin{eqnarray*}
|E'_r|=3k|E'_{r-1}|,~|V'_r|=|V'_{r-1}|+2k|E'_{r-1}|.
\end{eqnarray*}
It is easy to obtain
\begin{eqnarray*}
|E'_r|=m(3k)^r,~|V'_r|=n+2km\frac{(3k)^r-1}{3k-1}.
\end{eqnarray*}

 Denote $g_1(\lambda_i),$ $g_2(\lambda_i)$ and $g_3(\lambda_i)$ the roots of the equation $4\zeta^3-12\zeta^2+9\zeta-\lambda_i=0$. We at here define three new multiset as below.
\begin{eqnarray*}
g_1(U)=\bigcup_{x\in U}\{g_1(x)\},~~g_2(U)=\bigcup_{x\in U}\{g_2(x)\},~~g_3(U)=\bigcup_{x\in U}\{g_3(x)\}.
\end{eqnarray*}

In what follows, the normalized Laplacians of $S_{2k}^r(G)$ are given directly based on those notation and the structure of $S_{2k}^r(G)$.
\begin{thm} Assume that $S_{2k}^r(G)$ is the r-th iterations of the graph $S_{2k}(G)$, then the normalized Laplacian spectra of $S_{2k}^r(G)$ are as follows.
\begin{eqnarray*}
\begin{cases}
g_1(\Gamma(S_{2k}^{r-1}(G))\backslash \{0,2\})\bigcup g_2(\Gamma(S_{2k}^{r-1}(G))\backslash \{0,2\})\bigcup g_3(\Gamma(S_{2k}^{r-1}(G))\backslash \{0,2\})\\
\bigcup \{0,2\}\bigcup\bigg\{\underbrace{\frac{1}{2},\frac{1}{2}, \ldots, \frac{1}{2},\frac{1}{2}}_{k|E'_{r-1}|-|V'_{r-1}|+2}\bigg\}\bigcup\bigg\{\underbrace{\frac{3}{2},\frac{3}{2}, \ldots, \frac{3}{2},\frac{3}{2}}_{k|E'_{r-1}|-|V'_{r-1}|+2}\bigg\},~~~if~G~is~ bipartite,\\
g_1(\Gamma(S_{2k}^{r-1}(G))\backslash \{0\})\bigcup g_2(\Gamma(S_{2k}^{r-1}(G))\backslash \{0\})\bigcup g_3(\Gamma(S_{2k}^{r-1}(G))\backslash \{0\})\\
\bigcup \{0\}\bigcup\bigg\{\underbrace{\frac{1}{2},\frac{1}{2}, \ldots, \frac{1}{2},\frac{1}{2}}_{k|E'_{r-1}|-|V'_{r-1}|}\bigg\}\bigcup\bigg\{\underbrace{\frac{3}{2},\frac{3}{2}, \ldots, \frac{3}{2},\frac{3}{2}}_{k|E'_{r-1}|-|V'_{r-1}|+2}\bigg\},~~~if~G~is~nonbipartite.\\
\end{cases}
\end{eqnarray*}
\end{thm}

Based on the normalized Laplacians of $S_{2k}^r(G)$, the closed-form formulas for $Kf^*(S_{2k}^r(G))$, $Ke(S_{2k}^r(G))$ and $\tau(S_{2k}^r(G))$ are given.
\begin{thm}
Assume that $G$ is a $(n,m)$-graph. One has
\begin{eqnarray*}
Kf^*(S_{2k}^r(G))&=&(27k)^rKf^*(G)+\frac{16k^2(9k)^{r-1}(k^r-3^r)}{k-3}m^2+\frac{11k^2(3k)^{r-1}(9^r-1)}{8}m\\
&&-6k^2(3k)^{r-2}(9^r-1)n+\frac{4k^3(3k)^{r-2}[9-9^r+3k(9^r-1)-8(3k)^r]}{(3k-1)(k-3)}m^2.
\end{eqnarray*}
where $k$ is integer and $k\neq 3,~r\geq 1$.
\end{thm}
\noindent{\bf Proof.} According to Theorem 1.2 and $\lambda_i\neq 0,2$, one has $\zeta_1$, $\zeta_2$ and $\zeta_3$ are three roots of $4\zeta^3-12\zeta^2+9\zeta-\lambda_i=0$, then
\begin{eqnarray}
\frac{1}{\zeta_1}+\frac{1}{\zeta_2}+\frac{1}{\zeta_3}=\frac{9}{\lambda_i},~\zeta_1\zeta_2\zeta_3=\frac{\lambda_i}{4}.
\end{eqnarray}
In addition, $2+2+\frac{1}{2}=\frac{9}{2},~2\times \frac{1}{2}\times \frac{1}{2}=\frac{2}{4}=\frac{1}{2}$, hence $\lambda_i= 2$ also satisfy Eq.(5.22). Based on the definition of $Kf^*(G)$, one obtains
\begin{eqnarray}
Kf^*(S_{2k}^1(G))&=&6km\bigg(\sum_{i=2}^n\bigg(\frac{1}{\zeta_1}+\frac{1}{\zeta_2}+\frac{1}{\zeta_3}\bigg)+2(km-n)+\frac{2}{3}(km-n+2)
+\frac{1}{2}\bigg)\nonumber\\
&=&6km\bigg(\sum_{i=2}^n\frac{9}{\lambda_i}+2(km-n)+\frac{2}{3}(km-n+2)
+\frac{1}{2}\bigg)\nonumber\\
&=&27k\cdot Kf^*(G)+16k^2m^2-16kmn+11km.
\end{eqnarray}

By Eq.(5.23) and the construction of $S_{2k}^r(G)$, one has
\begin{eqnarray*}
Kf^*(S_{2k}^r(G))&=&27k\cdot Kf^*(S_{2k}^{r-1}(G))+16k^2|E'_{r-1}|^2-16k|E'_{r-1}||V'_{r-1}|+11k|E'_{r-1}|\\
&=&(27k)^rKf^*(G)+16k^2\sum_{i=0}^{r-1}(27k)^{r-1-i}|E'_i|^2\\
&&-16k\sum_{i=0}^{r-1}(27k)^{r-1-i}|E'_i||V'_i|
+11k\sum_{i=0}^{r-1}(27k)^{r-1-i}|E'_i|\\
&=&(27k)^rKf^*(G)+\frac{16k^2(9k)^{r-1}(k^r-3^r)}{k-3}m^2+\frac{11k^2(3k)^{r-1}(9^r-1)}{8}m\\
&&-6k^2(3k)^{r-2}(9^r-1)n+\frac{4k^3(3k)^{r-2}[9-9^r+3k(9^r-1)-8(3k)^r]}{(3k-1)(k-3)}m^2.
\end{eqnarray*}
This completes the proof.\hfill\rule{1ex}{1ex}

\noindent{\bf Remark 5.} When $k=3$, one has
\begin{eqnarray*}
Kf^*(S_6^1(G))=81\cdot Kf^*(G)+144m^2-48mn+33m.
\end{eqnarray*}
For $S_6^r(G)$, the multiplicative degree-Kirchhoff index is as follows.
\begin{eqnarray*}
Kf^*(S_{6}^r(G))&=&81\cdot Kf^*(S_{6}^{r-1}(G))+144|E'_{r-1}|^2-48|E'_{r-1}||V'_{r-1}|+33|E'_{r-1}|\\
&=&81^r\cdot Kf^*(G)+144\sum_{i=0}^{r-1}(81)^{r-1-i}|E'_i|^2-48\sum_{i=0}^{r-1}(81)^{r-1-i}|E'_i||V'_i|
+33\sum_{i=0}^{r-1}(81)^{r-1-i}|E'_i|\\
&=&81^r\cdot Kf^*(G)+\frac{11\cdot3^{2r-1}(9^r-1)}{8}m-2\cdot3^{2r-1}(9^r-1)mn\\
&&+\frac{3^{2r-1}[3(9^r-1)+8r\cdot9^r]}{2}m^2.
\end{eqnarray*}

\begin{thm}
Assume that $G$ is a $(n,m)$-graph. One has
\begin{eqnarray*}
Ke(S_{2k}^r(G))&=&9^rKe(G)+\frac{8k\cdot 3^{r-2}(k^r-3^r)}{k-3}m
-\frac{9^r-1}{3m}n+\frac{11k(9^r-1)}{48}\\
&&+\frac{2k[9-9^r+3k(9^r-1)-8(3k)^r]}{9(3k-1)(k-3)}m,
\end{eqnarray*}
where $k$ is integer and $k\neq 3,~r\geq 1$.
\end{thm}
\noindent{\bf Proof.} According to the relation between $Ke(G)$ and $Kf^*(G)$, one gets
\begin{eqnarray*}
Ke(S_{2k}^1(G))=\frac{1}{6km}Kf^*(S_{2k}^1(G))=9Ke(G)+\frac{1}{6}(16km-16n+11).
\end{eqnarray*}
Additionally,
\begin{eqnarray*}
Ke(S_{2k}^r(G))&=&\frac{1}{2m(3k)^r}Kf^*(S_{2k}^r(G))\\
&=&\frac{(27k)^r}{2m(3k)^r}Kf^*(G)+\frac{16k^2(9k)^{r-1}(k^r-3^r)}{2m(3k)^r(k-3)}m^2+\frac{11k^2(3k)^{r-1}(9^r-1)}{16m(3k)^r}m\\
&&-\frac{6k^2(3k)^{r-2}(9^r-1)}{2m(3k)^r}n+\frac{4k^3(3k)^{r-2}[9-9^r+3k(9^r-1)-8(3k)^r]}{2m(3k)^r(3k-1)(k-3)}m^2\\
&=&9^rKe(G)+\frac{8k\cdot 3^{r-2}(k^r-3^r)}{k-3}m-\frac{9^r-1}{3m}n+\frac{11k(9^r-1)}{48}\\
&&+\frac{2k[9-9^r+3k(9^r-1)-8(3k)^r]}{9(3k-1)(k-3)}m.
\end{eqnarray*}
As desired.\hfill\rule{1ex}{1ex}

\noindent{\bf Remark 6.} While $k=3$, this leads
\begin{eqnarray*}
Ke(S_6^r(G))&=&\frac{1}{2m\cdot9^r}Kf^*(S_6^r(G))\\
&=&9^rKe(G)+\frac{9^r-1}{48}(11-16n)+\frac{3(9^r-1)+8r\cdot9^r}{12}m.
\end{eqnarray*}

Before proceeding, we shall give an equation and define a function that will use in the following results. On the one hand, by Eq.(3.7), then
\begin{eqnarray}
\prod_{i=1}^{n+2km}d_i(S_{2k}(G))=2^{2km}\cdot\prod_{i=1}^{n}kd_i(G)=2^{2km}\cdot k^n\cdot\prod_{i=1}^{n}d_i(G).
\end{eqnarray}
On the other hand, we define
\begin{eqnarray*}
\varphi(r)=\frac{2km(r-1-3kr+3^rk^r)}{(3k-1)^2}.
\end{eqnarray*}

\begin{thm}
Assume that $G$ is a $(n,m)$-graph. One has
\begin{eqnarray*}
\tau(S_{2k}^r(G))=\bigg(\frac{1}{2}\bigg)^{r(1-n)-\varphi(r)}\bigg(\frac{3}{2}\bigg)^{r(2-n)-
\varphi(r)}3^{\frac{km(3^rk^r-1)}{3k-1}-1}k^{nr-1+\varphi(r)}\tau(G),
\end{eqnarray*}
where $k\geq 1,~r\geq 1$.
\end{thm}
\noindent{\bf Proof.}
By Eq.(5.24) and Lemma (2.4), one obtains
\begin{eqnarray*}
\tau(S_{2k}^1(G))&=&\frac{2\prod_{i=1}^{n+2km}d_i(S_{2k}(G))\prod^n_{i=2}\zeta_1\zeta_2\zeta_3\big
(\frac{1}{2}\big)^{km-n}\big(\frac{3}{2}\big)^{km-n+2}}{6km}\\
&=&\frac{2^{2km}\cdot k^n\cdot\prod_{i=1}^{n}d_i(G)\prod^n_{i=2}\lambda_i\big
(\frac{1}{2}\big)^{km-n+1}\big(\frac{3}{2}\big)^{km-n+2}}{6km}\\
&=&\bigg(\frac{1}{2}\bigg)^{1-n}\bigg(\frac{3}{2}\bigg)^{2-n}3^{km-1}k^{n-1}\cdot\tau(G).
\end{eqnarray*}
Moreover, the above equation yields
\begin{eqnarray*}
\tau(S_{2k}^r(G))&=&\bigg(\frac{1}{2}\bigg)^{1-|V'_{r-1}|}\bigg(\frac{3}{2}\bigg)^
{2-|V'_{r-1}|}3^{k|E'_{r-1}|-1}k^{|V'_{r-1}|-1}\cdot\tau(S_{2k}^{r-1}(G))\\
&=&\bigg(\frac{1}{2}\bigg)^{\sum_{i=0}^{r-1}(1-|V'_{i}|)}\bigg(\frac{3}{2}\bigg)^
{\sum_{i=0}^{r-1}(2-|V'_{i}|)}3^{k\sum_{i=0}^{r-1}(|E'_{i}|-1)}k^{\sum_{i=0}^{r-1}(|V'_{i}|-1)}\tau(G)\\
&=&\bigg(\frac{1}{2}\bigg)^{r(1-n)-\varphi(r)}\bigg(\frac{3}{2}\bigg)^{r(2-n)-\varphi(r)}3^
{\frac{km(3^rk^r-1)}{3k-1}-1}k^{nr-1+\varphi(r)}\tau(G).
\end{eqnarray*}
The result as desired. \hfill\rule{1ex}{1ex}

\section*{Acknowledgments}
%\ \ \ \ \textcolor{blue}{The authors would like to express their sincere gratitude to the
%editor and referees for many friendly
% and helpful suggestions, which led to great deal of improvements of the original
% manuscript.}

The work of was partly supported by the National Science Foundation of China under Grant Nos.
 11601006, 11801007, China Postdoctoral Science Foundation
Under Grant No. 2017M621579.

%, China Postdoctoral Science Foundation
%Under Grant No. 2017M621579, the Postdoctoral Science Foundation
%of Jiangsu Province Under Grant No. 1701081B, and Project of Anhui
%Jianzhu University under Grant 2016QD116 and Grant 2017dc03.

\end{document}